\documentclass[11pt,twoside,reqno,centertags,draft]{amsart}
\usepackage{amsfonts}
\usepackage{color,enumitem,graphicx}
\usepackage[colorlinks=true,urlcolor=blue,
citecolor=red,linkcolor=blue,linktocpage,pdfpagelabels,
bookmarksnumbered,bookmarksopen]{hyperref}

  \usepackage{amsmath,amsthm,amsfonts,amssymb}

\begin{document}

\title{  Anisotropic  singularities  to semi-linear elliptic equations in a measure framework  }
\date{}
\maketitle

\vspace{ -1\baselineskip}

{\small
\begin{center}

\bigskip

  {\sc  Huyuan Chen\footnote{chenhuyuan@yeah.net}}
  \bigskip

   Department of Mathematics, Jiangxi Normal University,

 \noindent Nanchang, Jiangxi 330022, PR China

\end{center}
}

\bigskip

\smallskip
\begin{quote}
{\bf Abstract.}
The purpose of this article is to study very  weak solutions of elliptic equation
\begin{equation}\label{00}
\arraycolsep=1pt\left\{
\begin{array}{lll}
 -\Delta    u+g(u)=2k\frac{\partial \delta_0}{\partial x_N }+j\delta_0\quad  &{\rm in}\quad\ \ B_1(0),\\[2mm]
 \phantom{-\Delta +g(u) }
 u=0\quad &{\rm on}\quad\ \  \partial B_1(0),
 \end{array}\right.
\end{equation}
where $k>0$, $j\ge0$,  $B_1(0)$ denotes the unit ball centered at the origin in $\mathbb{R}^N$ with $N\geq2$,  $g:\mathbb{R}\to\mathbb{R}$ is an odd, nondecreasing and $C^1$ function,
   $\delta_0$ is the Dirac mass concentrated at the origin
  and $\frac{\partial\delta_0}{\partial x_N}$ is defined in the distribution sense that
$$
\langle\frac{\partial \delta_0}{\partial x_N},\zeta\rangle=\frac{\partial\zeta(0)}{\partial x_N} , \qquad \forall \zeta\in C^1_0(B_1(0)).
$$
We obtain that problem (\ref{00}) admits a  unique very weak solution $u_{k,j}$ under  the integral subcritical assumption
$$\int_1^{\infty}g(s)s^{-1-\frac{N+1}{N-1}}ds<+\infty.$$
Furthermore, we prove that $u_{k,j}$  has anisotropic singularity
at the origin and we consider the odd property $u_{k,0}$ and limit of  $\{u_{k,0}\}_k$ as $k\to\infty$.

We pose the constraint on nonlienarity $g(u)$ that we only require integrability in the principle value sense, due to the singularities only  at the origin. This makes us able to search the very weak solutions in a larger scope of the nonlinearity.
\end{quote}

 \renewcommand{\thefootnote}{}
 \footnote{AMS Subject Classifications:  35R06, 35B40, 35Q60. }
\footnote{Key words:  Anisotropic singularity; Very weak solution; Uniqueness. }

\newcommand{\N}{\mathbb{N}}
\newcommand{\R}{\mathbb{R}}
\newcommand{\Z}{\mathbb{Z}}

\newcommand{\cA}{{\mathcal A}}
\newcommand{\cB}{{\mathcal B}}
\newcommand{\cC}{{\mathcal C}}
\newcommand{\cD}{{\mathcal D}}
\newcommand{\cE}{{\mathcal E}}
\newcommand{\cF}{{\mathcal F}}
\newcommand{\cG}{{\mathcal G}}
\newcommand{\cH}{{\mathcal H}}
\newcommand{\cI}{{\mathcal I}}
\newcommand{\cJ}{{\mathcal J}}
\newcommand{\cK}{{\mathcal K}}
\newcommand{\cL}{{\mathcal L}}
\newcommand{\cM}{{\mathcal M}}
\newcommand{\cN}{{\mathcal N}}
\newcommand{\cO}{{\mathcal O}}
\newcommand{\cP}{{\mathcal P}}
\newcommand{\cQ}{{\mathcal Q}}
\newcommand{\cR}{{\mathcal R}}
\newcommand{\cS}{{\mathcal S}}
\newcommand{\cT}{{\mathcal T}}
\newcommand{\cU}{{\mathcal U}}
\newcommand{\cV}{{\mathcal V}}
\newcommand{\cW}{{\mathcal W}}
\newcommand{\cX}{{\mathcal X}}
\newcommand{\cY}{{\mathcal Y}}
\newcommand{\cZ}{{\mathcal Z}}

\newcommand{\abs}[1]{\lvert#1\rvert}
\newcommand{\xabs}[1]{\left\lvert#1\right\rvert}
\newcommand{\norm}[1]{\lVert#1\rVert}

\newcommand{\loc}{\mathrm{loc}}
\newcommand{\p}{\partial}
\newcommand{\h}{\hskip 5mm}
\newcommand{\ti}{\widetilde}
\newcommand{\D}{\Delta}
\newcommand{\e}{\epsilon}
\newcommand{\bs}{\backslash}
\newcommand{\ep}{\emptyset}
\newcommand{\su}{\subset}
\newcommand{\ds}{\displaystyle}
\newcommand{\ld}{\lambda}
\newcommand{\vp}{\varphi}
\newcommand{\wpp}{W_0^{1,\ p}(\Omega)}
\newcommand{\ino}{\int_\Omega}
\newcommand{\bo}{\overline{\Omega}}
\newcommand{\ccc}{\cC_0^1(\bo)}
\newcommand{\iii}{\opint_{D_1}D_i}

\numberwithin{equation}{section}

\vskip 0.2cm \arraycolsep1.5pt
\newtheorem{lemma}{Lemma}[section]
\newtheorem{theorem}{Theorem}[section]
\newtheorem{definition}{Definition}[section]
\newtheorem{proposition}{Proposition}[section]
\newtheorem{remark}{Remark}[section]
\newtheorem{corollary}{Corollary}[section]

\setcounter{equation}{0}
\section{Introduction}
As early as in 1977, Lieb-Simon in \cite{LS} studied the very weak solutions to equation
\begin{equation}\label{eq 0.2}
 -\Delta u+(u-\lambda)_+^{\frac{3}{2}}=\sum_{i=1}^n m_i\delta_{a_i} \quad{\rm in}\quad \R^3
\end{equation}
in the description of  the Thomas-Fermi theory of electric field potential determined by the nuclear charge and distribution of electrons in an atom,
where $\lambda\ge0$, $t_+=\max\{t,0\}$,  $m_i>0$,  $a_i\in \R^3$ and $\delta_{a_i}$ is the Dirac mass at $a_i$ for $i=1,2,\cdots,n$. In fact, the solution of (\ref{eq 0.2}) turns out to be a
classical singular solution of
$$-\Delta u+(u-\lambda)_+^{\frac{3}{2}}=0 \quad{\rm in}\quad \R^3\setminus\{a_1,\cdots,a_n\}.$$
As a fundamental PDE's model,  the isolated singular problem
\begin{equation}\label{eq 0.1}
-\Delta u +|u|^{p-1} u= 0\quad {\rm in}\quad  \Omega\setminus \{0\}
\end{equation}
has been studied extensively, where $\Omega$ is a domain in $\R^N$ with $N\ge 3$.   Brezis-V\'eron in \cite{BV}
showed that problem (\ref{eq 0.1}) admits no isolated singular solution  when $p\ge \frac{N}{N-2}$.
A complete classification of the isolated singularities at the origin for (\ref{eq 0.1}) was given by  V\'{e}ron in \cite{V0} when $\frac{N+1}{N-1}\le p< \frac N{N-2}$ as follows:

 \noindent (i) either $|x|^{\frac{2}{p-1}}u(x)$  converges to a constant which can take only two values $\pm [\frac{2}{p-1}(\frac{2p}{p-1}-N)]^{1/(p-1)}$ as $x\to 0$,
  \smallskip

\noindent (ii) or $|x|^{N-2}u(x)$  converges to a constant $c_Nk$, and the couple $(u,\,k)$ is related to the weak solution of
$$
 -\Delta  u+u^p=k\delta_0 \quad   \rm{in}\quad (C^\infty_c(\Omega))',
$$
where $c_N$ is the normalized constant of the fundamental solution $\Gamma_N$ of  $-\Delta u=\delta_0$ in $\R^N$, that is,
\begin{equation}\label{Gamma}
 \Gamma_N(x)=\left\{
\arraycolsep=1pt
\begin{array}{lll}
 c_N|x|^{2-N} \quad  &{\rm if}\quad N\ge 3,\\[2mm]
 \phantom{  }
c_N \log|x|\quad &{\rm if}\quad N=2.
\end{array}
\right.
\end{equation}
 For  $1<p<\frac{N+1}{N-1}$,  the above classification holds under the restriction of nonnegative solutions of (\ref{eq 0.1})
 and  all above singular solutions  are isotropic. A conjucture states that there is a rich  structure of the  singularities for (\ref{eq 0.1})
 without the restriction of nonnegativity for $1<p<\frac{N+1}{N-1}$.
 V\'{e}ron in \cite{V0} partially answered this conjucture and showed that the anisotropic singular solutions could be constructed by considering the
 following nonlinear eigenvalue problem on $S^{N-1}$
 $$
  -\Delta_{S^{N-1}}\omega +|\omega|^{p-1}\omega =\lambda \omega,
$$
 where $S^{N-1}$ is the sphere of unit ball in $\R^N$ and $\Delta_{S^{N-1}}$ is the Laplace-Beltrami operator. Later on, Chen-Matano-V\'{e}ron in \cite{CMV} provideds the
   anisotropic singular solutions of (\ref{eq 0.1}) by analyzing the corresponding Laplace-Beltrami equations in the sphere. More singularities analysis see the references
   \cite{LS1,NS,RS0,V00}.

 In contrast with the absorption nonlinearity,  the isolated singular solutions of  elliptic  problem with source nonlinearity
\begin{equation}\label{eq 1.4}
 \left\{\begin{array}{lll}
 \displaystyle
  -\Delta   u=u^p\quad  {\rm in}\quad \Omega\setminus\{0\},
  \\[2mm]
 \displaystyle
 u> 0 \quad  {\rm in}\quad \Omega\setminus\{0\},\qquad
  u=0\quad
 {\rm on}\quad \partial\Omega
\end{array}\right.
\end{equation}
was classified  by Lions in \cite{L}, by using
Schwartz's Theorem to build that
\begin{equation}\label{3.3}
  -\Delta  u  -u^p =\sum_{|a|=0}^{\infty} k_a D^{a}\delta_0\quad{\rm in} \quad  ( C^\infty_c (\Omega))'
\end{equation}
 and then by choosing suitable test functions in $C^\infty_c (\Omega)$ to kill  all
 $D^{a}\delta_0$ with $a$ multiple index and $|a|\ge 1$, finally
building the connections with the weak solutions of
\begin{equation}\label{eq 1.5}
    \arraycolsep=1pt\left\{
\begin{array}{lll}
 \displaystyle    -\Delta    u=u^p+k\delta_0\quad
 &{\rm in}\quad \Omega,\\[2mm]
 \phantom{  -\Delta   }
 \displaystyle   u=0\quad
 &{\rm on}\quad  \partial\Omega.
\end{array}\right.
\end{equation}
Lions in \cite{L} proved that when $p\in (1,\frac{N}{N-2})$ with $N\ge3$, any  solution of (\ref{eq 1.4})
is a weak solution of (\ref{eq 1.5}) for some $k\ge0$, and when $p\ge \frac{N}{N-2}$, the parameter $k=0$.
Essentially, $D^{a}\delta_0$ with $|a|\ge 1$ is killed in (\ref{3.3}) because it is anisotropic singular source, that is, this source makes the solutions anisotropic singular.
For instance,  the fundamental solution of $-\Delta u=D^a\delta_0 $ with $a=(0,\cdots,0,1)\in\R^N$ is
\begin{equation}\label{1.4}
  P_N(x)= \tilde c_N\frac{x_N}{|x|^{N}},\quad \forall x\in\R^N\setminus\{0\},
\end{equation}
where   $\tilde c_N$ are the normalized constants, see \cite{CW}. Obviously, $P_N$ has anisotropic singularities.

Inspired by (\ref{3.3}), we observe that  $D^{a}\delta_0$ with $|a|\ge 1$ could provide anisotropic source and our motivation in this article  is to
make use of this kind of sources to construct  anisotropic singular solutions for (\ref{eq 0.1}) and to find the criteria for more general nonlinearity.
Let $B_1(0)$ be the unit ball centered at the origin  in $\R^N$ with $N\ge 2$,
denote by $\delta_0$ the Dirac
mass concentrated at the origin, to be convenient,  $\frac{\partial \delta_0}{\partial  x_N}=D^a\delta_0$ with $a=(0,\cdots,0,1)\in\R^N$
and then
$$\langle\frac{\partial \delta_0}{\partial  x_N},\zeta\rangle=\frac{\partial \zeta(0)}{\partial  x_N},\qquad \forall\, \zeta\in C^1_0(B_1(0)).$$
So our concern  is to study the isolated singular solutions of semilinear elliptic problem
\begin{equation}\label{eq 1.1}
\arraycolsep=1pt\left\{
\begin{array}{lll}
 -\Delta    u+g(u)=2k\frac{\partial \delta_0}{\partial x_N }+j\delta_0\quad  &{\rm in}\quad\ \ B_1(0),\\[2.5mm]
 \phantom{-\Delta +g(u) }
 u=0\quad &{\rm on}\quad\ \  \partial B_1(0),
\end{array}\right.
\end{equation}
where parameters $k>0$ and $j\ge0$.

Before we state our main results, we introduce the definition of   weak solution to (\ref{eq 1.1}) as follows.

\begin{definition}\label{weak sol}
A function $u\in
L^1(B_1(0))$ is called a very weak solution of (\ref{eq 1.1}),
 if  $g(u)$ be integrable in the principle value sense near the origin, $g(u)\in L^1(B_1(0), \, |x|dx)$ and 
\begin{equation}\label{1.10}
 \int_{B_1(0)} [u(-\Delta)\xi+ g(u)\xi]\, dx=2k\frac{\partial \xi(0)}{\partial  x_N}+j\xi(0),\qquad \forall\, \xi\in C^{1,1}_{0}(B_1(0)).
\end{equation}
\end{definition}

 We note that $k\frac{\partial \delta_0}{\partial x_N }$ and $j\delta_0$ are both visible in the distributional identity (\ref{1.10}).
 When $k=0$, the definition of very weak solution of (\ref{eq 1.1}) requires $g(u)\in L^1(B_1(0),\rho dx)$, see the references \cite{V}.  Since both sources have the support at the origin,
so  in this article we pose the constraint for the very weak solution   that $g(u)$ is integrable in the principle value sense near the origin means, i.e.  
$ \lim_{\epsilon\to0^+}\int_{B_1(0)\setminus B_\epsilon(0)} g(u)\, dx$ exists,
that  provides higher possibility for searching the sign-changing singular solutions of (\ref{eq 1.1}).

Now we are ready to state our first theorem on the existence and asymptotic behavior of very weak solutions to
problem (\ref{eq 1.1}).

\begin{theorem}\label{teo 1}
Assume that $k>0$, $j\ge0$, $\Gamma_N,$ $P_N$ is given in (\ref{Gamma}), (\ref{1.4}) respectively, and the nonlinearity $g:\R\to\R$ is an odd, nondecreasing and Lipschitz function satisfying
\begin{equation}\label{1.2}
\int_1^{\infty}g(s)s^{-1-\frac{N+1}{N-1}}ds<+\infty
\end{equation}
and
\begin{equation}\label{1.20}
g(s+t)-g(s)\le c_1 \left[\frac{g(s)}{1+|s|}t+g(t)\right],\quad s\in\R,\ t\ge0
\end{equation}
for some $c_1>0$.

Then (\ref{eq 1.1}) admits a  very weak solution $u_{k,j}$ such that

$(i)$ \ \ \ $u_{k,j}\ge 0$ in $ B_1^+(0)=\{x=(x',x_N)\in B_1(0):\, x_N>0\}$;

$(ii)$ \ \ $u_{k,j}$  has the following singularity at the origin
\begin{equation}\label{1.3}
\lim_{t\to0^+}\frac{u_{k,j}(te)}{P_{N}(te)}=2k\quad{\rm for}\ \ e=(e_1,\cdots,e_N)\in \partial B_1(0),\  e_N\not=0
\end{equation}
and
\begin{equation}\label{1.03}
\lim_{t\to0^+}\frac{u_{k,j}(te)}{\Gamma_{N}(te)}=j\quad{\rm for}\ \ e\in \partial B_1(0),\ e_N=0;
\end{equation}

$(iii)$ \ \ $u_{k,j}$ is a classical solution of
\begin{equation}\label{eq 1.2}
\arraycolsep=1pt\left\{
\begin{array}{lll}
 -\Delta    u+g(u)=0\quad  &{\rm in}\quad\ \ B_1(0)\setminus\{0\},\\[2mm]
 \phantom{-\Delta +g(u) }
 u=0\quad &{\rm on}\quad\ \  \partial B_1(0).
 \end{array}\right.
\end{equation}

In the particular case that $j=0$, denoting $u_k$ the solution $u_{k,0}$ of (\ref{eq 1.1}) with $j=0$, $u_k$ is $x_N$-odd, that is,
$$u_{k}(x',x_N)=-u_{k}(x',-x_N),\quad  \forall  (x',x_N)\in B_1(0)\setminus\{0\}.$$
Furthermore, the  $x_N$-odd very weak solution is unique.
\end{theorem}
\smallskip

We note that   $g(s)=|s|^{p-1}s$ with $p\in(0,\frac{N}{N-1})$ verifies (\ref{1.2}) and (\ref{1.20}). It follows by
(\ref{1.3}) that $|u_{k,j}|^{p-1}u_{k,j}\in L^1(B_1(0))$, but for $p\in[\frac{N}{N-1}, \frac{N+1}{N-1})$ and $k>0$,
$|u_{k,j}|^{p-1}u_{k,j}\not\in L^1(B_1(0))$. We can't able to obtain the uniqueness of the very weak solution to (\ref{eq 1.1}), due to the failure of application  the Kato's inequality, which requires that the nonlinearity term  $g(u)\in L^1(B_1(0))$.

  For the existence of very weak solutions, the normal method is to approximate  the Radon measure  by
$C^1_0$ functions and  consider the limit of the corresponding classical solutions.  When $k>0$,   we use  a sequence of Dirac measures $\frac{\delta_{te_N}-\delta_{-te_N}}{t}$ to approach  the source $\frac{\partial \delta_0}{\partial  x_N }$ and in this approximation, the biggest challenge is to
find a uniform estimate.  To overcome this difficulty, our strategy is to  to consider $x_N$-odd property of solutions when $j=0$ to derive the uniform bound in this approaching process.

We next state the nonexistence of very weak solution of (\ref{eq 1.1}).

\begin{theorem}\label{teo unique}
Assume that $k>0$, $j=0$, $g(s)=|s|^{p-1}s$ with $p\ge \frac{N+1}{N-1}$,
then there is no $x_N$-odd weak solution for problem (\ref{eq 1.1}).
\end{theorem}

Our strategy here is to  make use of $x_N-$odd property to deduce (\ref{eq 1.1}) into boundary data problem
 $$
\arraycolsep=1pt\left\{
\begin{array}{lll}
 -\Delta    u+g(u)=0\quad  &{\rm in}\quad\ \ B_1^+(0),\\[2mm]
 \phantom{-\Delta +g(u) }
 u=k\delta_0\quad &{\rm on}\quad\ \  \partial B_1^+(0),
\end{array}\right.
$$
in the distributional sense that
\begin{equation} \label{eq 1.6}
 \int_{B_1^+(0)} u(-\Delta)\xi dx+  \int_{B_1^+(0)} g(u)\xi dx=2k\frac{\partial \xi(x_0)}{\partial  x_N},\quad\quad\forall\, \xi\in C^{1,1}_0(B_1^+(0)).
\end{equation}
It is interesting but still open to derive the nonexistence when $j\not=0$.

Finally, we  analyze the limit of the weak solutions $\{u_{k}\}_k$ as $k\to\infty$.
From the monotonicity of $u_{k}$ in $B_1^+(0)$ and $B_1^-(0)$ respectively,  the limit of $\{u_{k}\}_k$ as $k\to\infty$ exists in $B_1(0)\setminus\{0\}$,
denote
\begin{equation}\label{1.01}
 u_\infty(x)=\lim_{k\to\infty} u_{k}(x),\qquad \forall\, x\in B_1(0)\setminus\{0\}.
\end{equation}

\begin{theorem}\label{teo 3}
 Assume that   $k>0$, $j=0$, $g(s)=|s|^{p-1}s$ with $p>1$, $u_{k}$ is the unique $x_N$-odd very weak solution of
(\ref{eq 1.1}) and $u_\infty$ is given by (\ref{1.01}).
Then    $u_\infty$ is a classical solution of
\begin{equation}\label{eq 1.03}
 \arraycolsep=1pt\left\{
\begin{array}{lll}
 -\Delta    u+|u|^{p-1}u=0\quad & {\rm in}\quad  B_1(0)\setminus\{0\},\\[2mm]
 \phantom{-\Delta    +|u|^{p-1}u}
u=0\quad & {\rm on}\quad \partial  B_1(0)
\end{array}\right.
\end{equation}
and satisfies that
\begin{equation}\label{1.5}
\lim_{t\to0^+} u_\infty(te)t^{\frac{2}{p-1}}= \varphi(e),\qquad \forall e\in \partial B_1(0),
\end{equation}
where
$\varphi:\partial B_1(0)\to \R$ is a continuous $x_N$-odd function such that for unit vector $e=(e_1,\cdots,e_N)$
$$\varphi(e)>0\quad {\rm if} \ \ e_N>0. $$

\end{theorem}

The rest of this paper is organized as follows. In Section 2, we analyze the $x_N$-odd property.   Section  3 is devoted to study the  $x_N$-odd very weak solution  in subcritical case when $j=0$ and the nonexistence  the  $x_N$-odd very weak solution  in the subcritical case.
In Section  4, we consider the limit of the unique $x_N$-odd weak solutions $u_{k}$ of (\ref{eq 1.1}) with $j=0$ as $k\to\infty$.
Finally, we prove the existence of non $x_N$-odd very weak solution when $j>0$ in Section 5.

\setcounter{equation}{0}
\section{Preliminary}

We start this section from the $x_N$-odd property.   Notice that   an  $x_N$-odd  function $w$ defined in a $x_N$-symmetric domain $B^*$ satisfies
$$
w(x)=0,\quad \forall x\in \{(x',0)\in B^*\}.
$$
In what follows,  we denote by $c_i$ a generic  positive constant.

\begin{lemma}\label{lm 2.1}
Assume that $f\in C^1(  \overline{B_1(0)})$ is an $x_N$-odd function, $g\in C^{1}(\R)$ is an odd  and nondecreasing
function.

 Then
 \begin{equation}\label{eq 2.1}
\arraycolsep=1pt\left\{
\begin{array}{lll}
 -\Delta    u+g(u)=f\quad  &{\rm in}\quad\ \ B_1(0),\\[2mm]
 \phantom{-\Delta +g(u) }
 u=0\quad &{\rm on}\quad\ \  \partial B_1(0)
 \end{array}\right.
\end{equation}
admits a unique classical solution $w_f$. Moreover,

$(i)$ $w_f$ is $x_N$-odd in $B_1(0)$;

$(ii)$  assume more that $f\ge0$ in $B_1^+(0)=\{x\in  B_1(0):\ x_N>0\}$ and $f\not\equiv0$ in $B_1^+(0)$,
then $w_f>0$ in $B_1^+(0)$.
\end{lemma}
{\bf Proof.} Since $g$  is an odd  and nondecreasing
function, then  it is standard to obtain the  existence of solution  by the method of super and sub solutions.

{\it Uniqueness.} Let $w_f,\tilde w_f$ be two solutions of (\ref{eq 2.1}),  $w=w_f-\tilde w_f$ in $B_1(0)$
and $A_+=\{x\in B_1(0):\, w(x)>0\}$. We claim that $A_+=\emptyset$. In fact, if $A_+\not=\emptyset$, we observe that  $w$
is a solution of
$$
\arraycolsep=1pt\left\{
\begin{array}{lll}
 -\Delta    w=g(\tilde w_f)-g( w_f)\le0\quad  &{\rm in}\quad\ \ A_+,\\[2mm]
 \phantom{-\Delta   }
 w=0\quad &{\rm on}\quad\ \  \partial A_+.
 \end{array}\right.
$$
By applying Maximum Principle, we have that
$$w\le 0\qquad{\rm in}\quad A_+,$$
which contradicts the definition of $A_+$. Then $A_+=\emptyset$.
Similarly, $\{x\in B_1(0): w(x)<0\} $ is empty.
Therefore, $w_f=\tilde w_f$ in $B_1(0)$ and the uniqueness holds.

{\it $(i)$ } Let $v(x',x_N)=-w_f(x',-x_N)$, and by direct computation, we derive that
\begin{eqnarray*}
 -\Delta v(x)+g(v(x))  &=& -\Delta[- w_f(x',-x_N)]+g(-w_f(x',-x_N)) \\
   &=&  \Delta w_f(x',-x_N)-g(w_f(x',-x_N))
   \\&=&-f(x',-x_N)=f(x),
\end{eqnarray*}
then $v$ is a solution of (\ref{eq 2.1}). It follows from the uniqueness of solution of (\ref{eq 2.1}) that
$$w_f(x',x_N)=-w_f(x',-x_N),\quad \forall x=(x',x_N)\in B_1(0).$$

{\it $(ii)$ }  We observe that
$w_f=0$ on $\partial B_1^+(0)$ and then
$w_f$ is a classical solution of
\begin{equation}\label{eq s52.3}
\arraycolsep=1pt\left\{
\begin{array}{lll}
 -\Delta    u+g(u)=f\quad  &{\rm in}\quad\ \ B_1^+(0),\\[2mm]
 \phantom{-\Delta +g(u) }
 u=0\quad &{\rm on}\quad\ \  \partial B_1^+(0).
 \end{array}\right.
\end{equation}
We now claim that $w_f\ge0$ in $ B_1^+(0)$. Indeed, if not, we have that $\min_{B_1^+(0)}w_f<0$.
Let $A_-=\{x\in B_1^+(0):\  w_f(x)<\frac12\min_{B_1^+(0)}w_f\}$,
  then $\psi:=w_f+\frac12\min_{B_1^+(0)}w_f$ satisfies
$$
\arraycolsep=1pt\left\{
\begin{array}{lll}
 -\Delta    \psi\ge 0\quad  &{\rm in}\quad\ \ A_-,\\[2mm]
 \phantom{-\Delta   }
\psi=0\quad &{\rm on}\quad\ \  \partial A_-.
 \end{array}\right.
 $$
By Maximum Principle, we have that
$$w_f(x)\ge \frac12\min_{B_1^+(0)}w_f,\quad x\in A_-,$$
which contradicts the definition of $A_-$.

We next prove that $w_f>0$ in $ B_1^+(0)$. Problem (\ref{eq s52.3}) could be seen as
$$\arraycolsep=1pt\left\{
\begin{array}{lll}
 -\Delta    u+\phi u=f\quad  &{\rm in}\quad\ \ B_1^+(0),\\[2mm]
 \phantom{-\Delta +gu }
 u=0\quad &{\rm on}\quad\ \  \partial B_1^+(0),
 \end{array}\right.
 $$
where $\phi(x)=\frac{g(w_f(x))}{w_f(x)}$ if $w_f(x)\not=0$ and $\phi(x)=g'(0)$ if $w_f(x)=0$. It follows by  $g\in C^1(\R)$  that $\phi$ is continuous and  $\phi\ge 0$ in $B_1^+(0)$.
Since $f\ge0$ in $B_1^+(0)$,  it follows by strong maximum principle  that $w_f>0$ or $w_f\equiv 0$ in $B_1^+(0)$, then we exclude $w_f\equiv 0$ in $B_1^+(0)$ by the fact that $f\not\equiv0$ in $B_1^+(0)$.
 \qquad$\Box$

\begin{corollary}\label{cr 2.1}
Assume that $f\in C( \overline{B_1(0)})$ is an $x_N$-odd function such that  $f\ge0$ in $B_1^+(0)$ and
 $g\in C^{1}(\R)$ is an odd  and nondecreasing function.
 Let $w_f$ be the solution of (\ref{eq 2.1}) and $\mathbb{G}_{B_1(0)}[f]$ be the
unique solution of
$$
\arraycolsep=1pt\left\{
\begin{array}{lll}
 -\Delta   u= f\quad  &{\rm in}\quad \ \ B_1(0),\\[2mm]
 \phantom{-\Delta   }
 u=0\quad &{\rm on}\quad  \partial B_1(0).
 \end{array}\right.
$$
Then $\mathbb{G}_{B_1(0)}[f]$ is $x_N$-odd and
$$0\le w_f\le \mathbb{G}_{B_1(0)}[f]\quad {\rm in}\quad B_1^+(0).$$

\end{corollary}
{\bf Proof.} By applying Lemma \ref{lm 2.1} with $g\equiv0$, we have that
$\mathbb{G}_{B_1(0)}[f]$ is $x_N$-odd and
$$\mathbb{G}_{B_1(0)}[f]\ge 0\quad {\rm in}\quad B_1^+(0).$$
Denote $v=\mathbb{G}_{B_1(0)}[f]-w_f$, then $v=0$ on $\partial B_1^+(0)$ and
$-\Delta v=g(w_f)\ge 0$, by Maximum Principle, we have that
$v\ge0$ in $B_1^+(0)$, which ends the proof. \qquad$\Box$

\begin{corollary}\label{cr 2.2}
Assume that $f\in C(  \overline{B_1(0)})$ is an $x_N$-odd function such that  $f\ge0$ in $B_1^+(0)$ and
 $ g_1,\, g_2\in C^{1}(\R)$   are odd  and nondecreasing functions satisfying
$$
 g_1(s)\le g_2(s),\quad \forall s\ge 0.
$$
Let $w_{f,i}$  be  the solutions of (\ref{eq 2.1}) replaced by $g$ by $g_i$ with $i=1,\,2$ respectively.

Then
$$|w_{f,1}(x)|\ge |w_{f,2}(x)|, \qquad \forall x\in B_1(0).$$

\end{corollary}
{\bf Proof.} By applying Lemma \ref{lm 2.1} and Corollary \ref{cr 2.1}, we have that
$w_{f,1}, w_{f,2}$ are $x_N$-odd and are nonnegative in $B_1^+(0)$.
We denote $w=w_{f,1}-w_{f,2}$, then $w$ satisfies that
$$
\arraycolsep=1pt\left\{
\begin{array}{lll}
 -\Delta w= g_2(w_{f,2})-g_1(w_{f,1}) \quad  &{\rm in}\quad\ \ B_1^+(0),\\[2mm]
 \phantom{-\Delta   }
 w=0\quad &{\rm on}\quad\ \  \partial B_1^+(0).
 \end{array}
 \right.
$$
We first claim that $w\ge0$ in $ B_1^+(0)$. If not, we have that $\min_{B_1^+(0)}w<0$.
Let us define $$A_-=\left\{x\in B_1^+(0):\ w(x)<\frac12\min_{B_1^+(0)}w\right\},$$
 then $\tilde w:=w+\frac12\min_{B_1^+(0)}w$ satisfies that
$$
\arraycolsep=1pt\left\{
\begin{array}{lll}
 -\Delta    \tilde w\ge 0\quad  &{\rm in}\quad\ \ A_-,\\[2mm]
 \phantom{-\Delta   }
\tilde  w=0\quad &{\rm on}\quad\ \  \partial A_-.
 \end{array}\right.
 $$
By Maximum Principle, we have that
$$  w(x)\ge \frac12\min_{B_1^+(0)}w,\quad \forall x\in A_-,$$
which contradicts the definition of $A_-$. \qquad$\Box$

\begin{proposition}\label{pr 2.1}

 Let   $f_1$ and $ f_2$ be  $x_N$-odd functions in $ C^1_{loc}( \overline{B_1(0)}\setminus\{0\})\cap L^1(B_1(0),|x| dx)$ satisfying $f_2\ge f_1\ge 0$ in $B_1^+(0)$,  then
the problem
\begin{equation}\label{homo}
 \arraycolsep=1pt\left\{
\begin{array}{lll}
 -\Delta  &u=f_i\quad & {\rm in}\quad B_1(0),\\[2mm]
&u=0\quad & {\rm on}\quad \partial B_1(0)
\end{array}\right.
\end{equation}
admits a unique  $x_N$-odd weak solution $u_i$ with $i=1,2$ in the sense that $u_i\in L^1(B_1(0))$,
\begin{equation}\label{2.6}
 \int_{B_1(0)} u_i(-\Delta)\xi  dx= \int_{B_1(0) } \xi f_idx ,\quad \forall \xi\in C^{1,1}_{0}(B_1(0)),\ \xi(0)=0.
\end{equation}
Moreover, $u_i$ is a classical solution of
\begin{equation}\label{v 1}
\arraycolsep=1pt\left\{
\begin{array}{lll}
 -\Delta  &u=f_i \quad & {\rm in}\quad B_1(0)\setminus\{0\},\\[2mm]
&u=0\quad & {\rm on}\quad \partial B_1(0)
\end{array}\right.
\end{equation}
and
 \begin{equation}\label{sign}
0\le u_1(x)\le u_2(x)\le  \int_{B_1^+(0)}  \frac{c_2 f_2(y)|y|}{|y-x| (|y-x|+2|y|)^{N-2}}dy,\quad \forall x\in B_1^+(0).
\end{equation}

\end{proposition}
{\bf Proof.} {\it Uniqueness}.
Let $w$ satisfy that
\begin{equation}\label{L0}
 \int_{B_1(0)} w(-\Delta)\xi  dx= 0,
\end{equation}
for any $\xi\in C^{1,1}_{0}(B_1(0))$ such that $\xi(0)=0$. Since $0\in L^1(B_1(0))$, then
the test function could be improved into $C^{1,1}_0(B_1(0))$ without the restriction that $\xi(0)=0$.
Denote by $\eta_1$ the
solution of
\begin{equation}\label{2.8}
\left\{
\begin{array}{lll}
-\Delta \eta_1 ={\rm sign}(w)&\quad\mbox{ in }\quad B_1(0),\\[2mm]
\phantom{-\Delta }
\eta_1=0&\quad\mbox{ on
}\quad\partial B_1(0).
\end{array}\right.
\end{equation}
Then $\eta_1\in C^{1,1}_{0}(B_1(0))$ and  then
$$\displaystyle
\int_{B_1(0)}|w|  \ dx=0.
$$
 This implies $w=0$ in $B_1(0)\setminus\{0\}$.
\smallskip
\smallskip

{\it Existence}. Let $f_{i,\epsilon}=f_i \chi_{B_1(0)\setminus B_\epsilon(0)}$,  where $i=1,2$, $\chi_{B_1(0)\setminus B_\epsilon(0)}=1$ in $B_1(0)\setminus B_\epsilon(0)$ and  $\chi_{B_1(0)\setminus B_\epsilon(0)}=0$ in $B_\epsilon(0)$,  then  $f_{i,\epsilon}$ is an $x_N$-odd function in  $L^{\infty}(B_1(0))$ such that $f_{2,\epsilon}\ge f_{1,\epsilon}\ge 0$
 in $B^+_1(0)$, then
\begin{equation}\label{homo 1}
 \arraycolsep=1pt\left\{
\begin{array}{lll}
 -\Delta  &u=f_{i,\epsilon}\quad & {\rm in}\quad B_1(0),\\[2mm]
&u=0\quad & {\rm on}\quad \partial B_1(0)
\end{array}\right.
\end{equation}
admits a unique solution $u_{i,\epsilon}$
satisfying
\begin{equation}\label{2.007}
 \int_{B_1(0)} u_{i,\epsilon}(-\Delta)\xi  dx= \int_{B_1(0) } \xi f_{i,\epsilon} dx ,\qquad \forall \xi\in C^{1,1}_{0}(B_1(0)),\ \xi(0)=0.
\end{equation}
Moreover, from Lemma \ref{lm 2.1} with $g\equiv0$, we have that
$$0\le u_{i,\epsilon} \le u_{i,\epsilon'}\quad {\rm in} \ B_1^+(0) \quad {\rm for }\quad 0<\epsilon'\le \epsilon, $$
$$0\le u_{1,\epsilon} \le u_{2,\epsilon}\quad {\rm in} \ B_1^+(0) \quad {\rm for \ any}\quad  \epsilon\ge 0, $$
and for any $x\in B_1^+(0)$,
\begin{eqnarray*}
  u_{2,\epsilon}(x)& =& \int_{B_1(0)}G_{B_1(0)}(x,y)  f_{2,\epsilon}(y) dy \\
   &=&\int_{B_1^+(0)}[G_{B_1(0)}(x,y)-G_{B_1(0)}(x,\tilde y)]   f_{2,\epsilon}(y)dy,
\end{eqnarray*}
where $\tilde y=(y',-y_N)$.
We observe that $G_{B_1(0)}(x,y)=\frac{c_N}{|x-y|^{N-2}}-\tilde G_{B_1(0)}(x,y)$, where $\tilde G_{B_1(0)}(x,y)$ is a harmonic function in $B_1(0)$ with
the boundary value $\frac{c_N}{|x-y|^{N-2}}$ for $y\in \partial B_1(0)$. Therefore, for $x,y\in B_1^+$, we have that $\tilde G_{B_1(0)}(x,y)-\tilde G_{B_1(0)}(x,\tilde y)\ge0$
and then
$$G_{B_1(0)}(x,y)-G_{B_1(0)}(x,\tilde y)\le c_3\left[\frac1{|y-x|^{N-2}}-\frac1{|\tilde y-x|^{N-2}}\right].$$
 Moreover,   we see that
$$|\tilde y-x|\le |y-x|+2y_N\le|y-x|+2|y| $$
and
\begin{eqnarray*}
 0 \le   \frac1{|y-x|^{N-2}}-\frac1{|\tilde y-x|^{N-2}}
   &\le & \frac1{|y-x|^{N-2}}-\frac1{(|y-x| +2|y|)^{N-2}}
 \\ &=&  \frac{(|y-x|+2|y|)^{N-2}-|y-x|^{N-2}}{|y-x|^{N-2}(|y-x|+2|y|)^{N-2}}
 \\&\le & \frac{c_4|y|}{|y-x| (|y-x|+2|y|)^{N-2}},
\end{eqnarray*}
then
\begin{eqnarray*}
u_2(x)& \le & c_2\int_{B_1^+(0)}  \frac{|y|f_2(y)}{|y-x| (|y-x| +2|y|)^{N-2}}dy.
\end{eqnarray*}
Therefore, we obtain a uniform bound for $u_{2,\epsilon}$, together with monotonicity, then passing to the limit as $\epsilon \to0^+$ in (\ref{2.007}), we deduces that
$u_i:=\lim_{\epsilon\to0^+}u_{i,\epsilon}$ is a weak solution of (\ref{homo}) and it follows from standard stability theorem that
$ u_{i}$ is a classical solution of  (\ref{v 1}) for $i=1,2$.\quad $\Box$

\smallskip

By direct extension, we have the  following corollary.

\begin{corollary}\label{pr 2.2}

 Let   $f$   be an $x_N$-odd function in $ C^1_{loc}( \overline{B_1(0)}\setminus\{0\})\cap L^1(B_1(0),|x|^i dx)$ with $i\in\Z$ satisfying $f \ge  0$ in $B_1^+(0)$,  then
the problem
$$
 \arraycolsep=1pt\left\{
\begin{array}{lll}
 -\Delta  &u=f\quad & {\rm in}\quad B_1(0),\\[2mm]
&u=0\quad & {\rm on}\quad \partial B_1(0)
\end{array}\right.
$$
admits a unique  $x_N$-odd weak solution $u_f$ with $i\in\N$, $i\ge 2$ in the sense that $u_f\in L^1(B_1(0),|x|^{i-1}dx)$,
$$
 \int_{B_1(0)} u(-\Delta)\xi  dx= \int_{B_1(0) } \xi fdx,
$$
for any $\xi\in C^{1,1}_{0}(B_1(0))$ s. t.  $|\xi(x)|\le c|x|^{i}$ for any $x\in B_1(0)$ and some $c>0$.
Moreover,
$$
 u_f(x)\ge 0,\quad \forall x\in B_1^+(0).
$$
and $u_f$ is a classical solution of
$$
\arraycolsep=1pt\left\{
\begin{array}{lll}
 -\Delta  &u=f \quad & {\rm in}\quad B_1(0)\setminus\{0\},\\[2mm]
&u=0\quad & {\rm on}\quad \partial B_1(0).
\end{array}\right.
$$

\end{corollary}

\begin{remark}
The arguments in Proposition \ref{pr 2.1} and Corollary  \ref{pr 2.2} hold when  $B_1(0)$ is replaced by $\R^N$ and
the boundary condition is done  by
$$\lim_{|x|\to+\infty} u(x)=0.$$

\end{remark}

In order to study the convergence of weak solutions, we recall the definition and basic properties of the Marcinkiewicz
spaces.

\begin{definition}\label{definition 2.1}
Let $\Omega\subset \R^N$ be a domain and $\mu$ be a positive
Borel measure in $\Omega$. For $\kappa>1$,
$\kappa'=\kappa/(\kappa-1)$ and $u\in L^1_{loc}(\Omega,d\mu)$, we
set
$$
 \arraycolsep=1pt
\begin{array}{lll}
\|u\|_{M^\kappa(\Omega,d\mu)}=\inf\{c\in[0,\infty]:\int_E|u|d\mu\le
c\left(\int_Ed\mu\right)^{\frac1{\kappa'}},\ \forall E\subset\Omega\
{\rm Borel\ set}\}
\end{array}
$$
and
\begin{equation}\label{spa M}
M^\kappa(\Omega,d\mu)=\{u\in
L_{loc}^1(\Omega,d\mu):\|u\|_{M^\kappa(\Omega,d\mu)}<\infty\}.
\end{equation}
\end{definition}
$M^\kappa(\Omega,d\mu)$ is called the Marcinkiewicz space with
exponent $\kappa$ or weak $L^\kappa$ space and
$\|.\|_{M^\kappa(\Omega,d\mu)}$ is a quasi-norm. The following
property holds.

\begin{proposition}\label{pr 1} \cite{BBC}
Assume that $1\le q< \kappa<\infty$ and $u\in L^1_{loc}(\Omega,d\mu)$.
Then there exists  $C(q,\kappa)>0$ such that
$$\int_E |u|^q d\mu\le C(q,\kappa)\|u\|_{M^\kappa(\Omega,d\mu)}\left(\int_E d\mu\right)^{1-q/\kappa},$$
for any Borel set $E$ of $\Omega$.
\end{proposition}

\setcounter{equation}{0}
\section{$x_N$-odd very weak solution with $j=0$ }

\subsection{Existence of very weak solution}

In this subsection, we prove the existence and uniqueness of very weak solution to problem (\ref{eq 1.1}) when $j=0$.

\begin{theorem}\label{teo 3.1}
Assume that $k>0$ $j=0$,   the nonlinearity $g:\R\to\R$ is an odd, nondecreasing and Lipchitz continuous function satisfying
(\ref{1.2}).

Then (\ref{eq 1.1}) admits a unique  $x_N$-odd very weak solution $w_{k}$ such that

$(i)$ \ \ \ $w_{k'}\ge w_{k}\ge 0$ in $ B_1^+(0)$ for $k'\ge k>0$;

$(ii)$ \ \ $w_{k}$  satisfies (\ref{1.3});

$(iii)$ \ \ $w_{k}$ is a classical solution of (\ref{eq 1.2}).

\end{theorem}

Before proving Theorem \ref{teo 3.1}, we need following preliminaries.

\begin{lemma}\label{lm 3.1}
Assume that $k>0$,  the nonlinearity $g:\R\to\R$ is an odd, nondecreasing and Lipchitz continuous function and $u$ is a very weak solution of (\ref{eq 1.1}), locally bounded in $B_1(0)\setminus\{0\}$.
Then $u$ is a classical solution of (\ref{eq 1.2}).

\end{lemma}
{\bf Proof.} Since $\frac{\partial \delta_0}{\partial x_N },\, \delta_0$ have the support in $\{0\}$, so for any open sets $O_1,O_2$ in $B_1(0)$ such that
 $\bar O_1\subset O_2\subset \bar O_2\subset B_1(0)\setminus\{0\},$
$u$ is
a very weak solution of
\begin{equation}\label{2.09}
 -\Delta u+g(u)=0\quad {\rm in}\quad O_2,
\end{equation}
where $u\in L^\infty(O_2)$ and $g(u)\in L^\infty(O_2)$.
By standard regularity results, we have that
$u$ satisfies (\ref{2.09}) in $O_1$ in the classical sense. \qquad$\Box$

\smallskip
For $n\in\N$, we consider $\{g_n\}_n$ of $C^1$, odd, nondecreasing  functions defined in $\R$
satisfying
\begin{equation}\label{3.1}
  g_n\le g_{n+1}\le g\ {\rm in}\ \R_+,\quad \sup_{s\in\R_+}g_n(s)=n\quad{\rm and}\quad \lim_{n\to\infty}\norm{g_n-g}_{L^\infty_{loc}(\R)}=0.
\end{equation}

\begin{proposition}\label{pr 3.1}
Let $g_n$ be defined by (\ref{3.1}) and
$$\mu_t=\frac{\delta_{te_N}-\delta_{-te_N}}{t},\quad \ t\in(0,1).$$
Then for any $n\in\N$, problem
\begin{equation}\label{eq s52.6}
\arraycolsep=1pt\left\{
\begin{array}{lll}
 -\Delta    u+g_n(u)=k\mu_t\quad  &{\rm in}\quad\ \ B_1(0),\\[2mm]
 \phantom{-\Delta +g_n(u) }
 u=0\quad &{\rm on}\quad\ \  \partial B_1(0)
 \end{array}\right.
\end{equation}
admits a unique very weak solution $w_{k,n,t}$, which is a classical solution of
$$
\arraycolsep=1pt\left\{
\begin{array}{lll}
 -\Delta    u+g_n(u)=0\quad  &{\rm in}\quad\ \ B_1(0)\setminus\{te_N,-te_N\},\\[2mm]
 \phantom{-\Delta +g_n(u) }
 u=0\quad &{\rm on}\quad\ \  \partial B_1(0).
 \end{array}\right.
$$
 Moreover,   $w_{k,n,t}$ is $x_N$-odd in $B_1(0)\setminus\{te_N,-te_N\}$,
$$
0\le w_{k,n+1,t}\le w_{k,n,t} \le \mathbb{G}_{B_1(0)} [\mu_t]\quad {\rm in}\quad B_1^+(0)\setminus\{te_N\}
$$
and
\begin{equation}\label{s52.8}
 w_{k+1,n,t}\ge w_{k,n,t}  \quad {\rm in}\quad B_1^+(0)\setminus\{te_N\}.
\end{equation}
\end{proposition}
{\bf Proof.}  We observe that $\mu_t$ is a bounded Radon measure and $g_n$ is bounded, Lipschitz continuous and nondecreasing, then it follows from \cite[Theorem  3.7]{V} under the integral subcritical
assumption (\ref{1.2}) replaced $\frac{N}{N-1}$ by $\frac{N}{N-2}$ for $N\ge3$ and the Kato's inequality, problem (\ref{eq s52.6}) admits a unique weak solution $w_{k,n,t}$. Moreover, $w_{k,n,t}$ could be approximated by
the classical solutions $\{w_{k,n,t,m}\}$ to problem
\begin{equation}\label{3.2}
\arraycolsep=1pt\left\{
\begin{array}{lll}
 -\Delta    u+g_n(u)=k\mu_{t,m}\quad  &{\rm in}\quad\ \ B_1(0),\\[2mm]
 \phantom{-\Delta +g_n(u) }
 u=0\quad &{\rm on}\quad\ \  \partial B_1(0),
 \end{array}\right.
\end{equation}
where $$\mu_{t,m}(x)=\frac{\sigma_m(x-te_N)-\sigma_m(x+te_N)}t$$ and
$\{\sigma_m\}_m$ is a sequence of radially symmetric, nondecreasing smooth functions converging
to $\delta_0$ in the distribution sense. Furthermore,
\begin{equation}\label{s52.01}
\int_{B_1(0)} [w_{k,n,t,m}(-\Delta)\xi+ g_n(w_{k,n,t,m})\xi] dx=k\int_{B_1(0)}\mu_{t,m}\xi dx,\quad \forall \xi\in C^{1,1}_0(B_1(0)).
\end{equation}

Since $\mu_{t,m}$ is $x_N$-odd and  nonnegative in $B_1^+(0)$,  so is $w_{k,n,t,m}$ by Lemma \ref{lm 2.1}.
We observe that $(k+1)\sigma_m\ge k \sigma_m$, it follows from Lemma \ref{lm 2.1} that
$$
w_{k+1,n,t,m}\ge w_{k,n,t,m}\quad {\rm in}\ B_1^+(0).
$$
Since $$g_n(s)s\le g_{n+1}(s)s,\quad \forall s\in\R,$$
then it follows from Corollary \ref{cr 2.2}  that
\begin{equation}\label{s52.02}
w_{k,n,t,m}\le w_{k,n+1,t,m}\quad {\rm in}\ B_1^+(0).
\end{equation}
From the proof of Theorem 2.9 in \cite{V}, we know  that
$$\mathbb{G}_{B_1(0)}[\mu_{t,m}]\to \mathbb{G}_{B_1(0)}[\mu_{t}]\quad {\rm in}\quad B_1(0)\setminus\{te_N,-te_N\}\quad {\rm and\quad in}\ L^q({B_1(0)})\quad {\rm as }\  m\to\infty,$$
where $q\in[1,\frac{N}{N-2})$.
By regularity results, any compact set $K$ and open set $O$ in $B_1(0)$ such that $K\subset O$, $\bar O\cap \{te_N,-te_N\}=\emptyset$, there exist $c_6,c_7>0$ independent of $m$ such that
\begin{eqnarray*}
 \norm{w_{k,n,t,m}}_{C^2(K)} &\le & c_6\norm{k\mathbb{G}_{B_1(0)}[\mu_{t,m}] }_{L^\infty(O)} \\
   &\le & c_7k\norm{\mathbb{G}_{B_1(0)}[\mu_{t}] }_{L^\infty(O)}.
\end{eqnarray*}
Therefore, up to some subsequence, there exists a measurable function  $\tilde w$ such that
$$w_{k,n,t,m}\to \tilde w\quad {\rm in}\quad B_1(0)\setminus\{te_N,-te_N\}\quad {\rm and\quad in}\ L^q(B_1(0))\quad {\rm as }\  m\to\infty,$$
where $q\in[1,\frac{N}{N-2})$.
Then $\tilde w$ is $x_N$-odd, nonnegative in $B_1^+(0)$ and
$$g_n(w_{k,n,t,m})\to g_n(\tilde w)\quad {\rm in}\quad B_1(0)\setminus\{te_N,-te_N\}\quad {\rm and\quad in}\ L^1(B_1(0))\quad {\rm as }\  m\to\infty.$$
Passing to the limit in (\ref{s52.01}) as $m\to\infty$, we deduce that $\tilde w$ is
a weak solution of (\ref{eq s52.6}).
By the uniqueness of weak solution of (\ref{eq s52.6}), we obtain that $w_{n,t}=\tilde w$.
Therefore, $w_{n,t}$ is $x_N$-odd, nonnegative in $B_1^+(0)$ and it follows from (\ref{s52.02}) and (\ref{s52.8}) that
$$
w_{k,n,t}\le w_{k,n+1,t}\quad {\rm in}\ \ B_1^+(0)
$$
and
$$
w_{k+1,n,t}\ge w_{k,n,t}\quad {\rm in}\ \ B_1^+(0).
$$
This ends the proof.\qquad $\Box$

\smallskip

We next passing to the limit of   weak solutions as $t\to0^+$.

\begin{proposition}\label{pr s52.2}
Let $g_n$ be defined by (\ref{3.1}).
Then for any $n\in\N$, problem
\begin{equation}\label{eq 3.1}
\arraycolsep=1pt\left\{
\begin{array}{lll}
 -\Delta    u+g_n(u)=2k\frac{\partial \delta_0}{\partial x_N}\quad  &{\rm in}\quad\ \ B_1(0),\\[2mm]
 \phantom{-\Delta +g_n(u) }
 u=0\quad &{\rm on}\quad\ \  \partial B_1(0)
 \end{array}\right.
\end{equation}
admits a unique very weak solution $w_{k,n}$. Moreover, \\
$(i)$  $w_{k,n}$ is $x_N$-odd for any $n\in\N$ in $B_1(0)\setminus\{0\}$ and
$$
0\le w_{k,n+1}\le w_{k,n} \le 2k\mathbb{G}_{B_1(0)} [\frac{\partial \delta_0}{\partial x_N}]\quad {\rm in}\quad B_1^+(0)
$$
and
$$
0\le w_{k,n}\le w_{k+1,n} \quad {\rm in}\quad B_1^+(0);
$$
$(ii)$ $w_{k,n}$ is a classical solution of
$$
\arraycolsep=1pt\left\{
\begin{array}{lll}
 -\Delta    u+g_n(u)=0\quad  &{\rm in}\quad\ \ B_1(0)\setminus\{0\},\\[2mm]
 \phantom{-\Delta +g_n(u) }
 u=0\quad &{\rm on}\quad\ \  \partial B_1(0).
 \end{array}\right.
$$
\end{proposition}
{\bf Proof.}  It follows from Proposition \ref{pr 3.1} that problem (\ref{eq s52.6}) admits a unique very weak solution $w_{k,n,t}$,
that is,
\begin{equation}\label{q s52.2}
\int_{B_1(0)} [w_{k,n,t}(-\Delta)\xi+g_n(w_{k,n,t}) ]dx=k\frac{\xi(te_N)-\xi(-te_N)}{t},\quad \forall \xi\in  C^{1,1}_0(B_1(0)).
\end{equation}
On the one hand,  we have that
$$\lim_{t\to0^+}\frac{\xi(te_N)-\xi(-te_N)}{t}=2\frac{\partial\xi(0)}{\partial x_N}.$$
On the other hand, by Proposition \ref{pr 3.1}, we have that
$$|w_{k,n,t}|\le k|\mathbb{G}_{B_1(0)}[\mu_t]|\quad{\rm in}\ B_1(0).$$
By regularity results, for $\sigma\in(0,1)$ and any compact set $K$ and open set $O$ in $B_1(0)$ such that $K\subset O$, $\bar O\cap \{te_N:\ t\in(-\frac12,\frac12)\}=\emptyset$,
there exist $c_8,c_9>0$ independent of $t$ such that
\begin{eqnarray*}
 \norm{w_{k,n,t}}_{C^{2+\sigma}(K)} \le  c_8\norm{k\mathbb{G}_{B_1(0)}[\mu_{t}] }_{L^\infty(O)}
   \le  c_9k\norm{\mathbb{G}_{B_1(0)}[\frac{\partial\delta_{0}}{\partial x_N}]}_{L^\infty(O)}.
\end{eqnarray*}
Moreover, by \cite[Proposition 3.3]{CW}, $\{\mathbb{G}_{B_1(0)}[\mu_t]\}$ is uniformly bounded in $M^{\frac{N}{N-1}}(B_1(0),dx)$ if $N\ge 3$
 and is uniformly bounded in $M^{\frac{2}{1+\sigma}}(B_1(0),dx)$ for any $\sigma\in(0,\frac12)$ if $N=2$.
 Therefore, $\{w_{k,n,t} \}_t$ is relatively compact in $L^p({B_1(0)})$ for any $p\in[1,\frac{N}{N-1})$.
There exists $w_{k,n}\in L^1(B_1(0))$ such that
$$w_{k,n,t}\to w_{k,n}\quad {\rm  a.e.\ in}\ B_1(0)\quad {\rm and\  in}\  L^1(B_1(0)),$$
which implies that
$$g_n(w_{k,n,t})\to g_n(w_{k,n})\quad {\rm  a.e.\ in}\ B_1(0)\quad {\rm and\  in}\  L^1(B_1(0))\ {\rm as}\ t\to0^+.$$

 Therefore,  up
to some subsequence, passing to the limit as $t\to0^+$ in the identity (\ref{q s52.2}), it follows that
$w_{k,n}$ is a very  weak solution of (\ref{eq 3.1}). Moreover, $w_{k,n}$ is $x_N$-odd and nonnegative in $B_1^+(0)$.

{\it Uniqueness.} Let $v_n$ be a weak solution of (\ref{eq 3.1}) and then $\varphi_n:=w_{k,n}-v_n$ is a very weak solution to
$$
\arraycolsep=1pt\left\{
\begin{array}{lll}
-\Delta \varphi_n +g_n(w_{k,n})-g_n(v_n)=0\quad &{\rm in}\quad B_1(0),
\\[2mm]\phantom{-\Delta_n +g(w_{k,n})-g_n(v_n)}
\varphi_n=0\quad &{\rm on}\quad \partial B_1(0).
\end{array}\right.
$$
By Kato's inequality \cite[Theorem 2.4]{V} (see also \cite{K,LL0,RS0}),
$$\int_{B_1(0)}|\varphi_n|(-\Delta)\xi+\int_{B_1(0)}[g_n(w_{k,n})-g(v_n)]{\rm sign}(w_{k,n}-v_n)\xi\, dx\le 0 $$
Taking $\xi=\mathbb{G}_{B_1(0)}[1]$, we have that
$$\int_{B_1(0)}[g_n(w_{k,n})-g(v_n)]{\rm sign}(w_{k,n}-v_n)\xi\, dx\ge 0\quad{\rm and}\quad \int_{B_1(0)}|\varphi_n|\,dx=0,$$
then $\varphi_n=0$ a.e. in $B_1(0)$. Then the uniqueness is obtained. \qquad$\Box$

\smallskip

The next estimate plays an important role in $w_{k,n}\to w_k$ in $L^p({B_1(0)})$ with $p\in[1,\frac{N+1}{N-1})$.

\begin{lemma}\label{lm s53.1}

 There exists $c_{10}>0$ such that
 \begin{equation}\label{s53.0}
\| \mathbb{G}_{B_1(0)}[\frac{\partial\delta_{0}}{\partial x_N}]\|_{M^{\frac{N}{N-1}}(B_1(0))}\le c_{10}
\end{equation}
and
\begin{equation}\label{s53.1}
\| \mathbb{G}_{B_1(0)}[\frac{\partial\delta_{0}}{\partial x_N}]\|_{M^{\frac{N+1}{N-1}}(B_1(0), |x|dx)}\le c_{10}.
\end{equation}

\end{lemma}
{\bf Proof.} We  observe that
 $$\mathbb{G}_{B_1(0)}[\frac{\partial\delta_{0}}{\partial x_N}](x)= \frac{\partial G_{B_1(0)}(x,0)}{\partial x_N}$$
and  for $x,y\in{B_1(0)}$, $x\not=y$,
$$
G_{B_1(0)}(x,y)=\left\{\arraycolsep=1pt
\begin{array}{lll}
c_N|x-y|^{2-N}+\tilde G_{B_1(0)}(x,y)\quad &{\rm if}\quad N\ge 3,
\\[2mm]
-c_N\log|x-y|+\tilde G_{B_1(0)}(x,y)\quad&{\rm if}\quad N=2,
\end{array}
\right.
$$
where $\tilde G_{B_1(0)}$ is a harmonic function in $B_1(0)\times B_1(0)$.
Then
\begin{eqnarray*}
  \left|\frac{\partial G_{B_1(0)}(x,0)}{\partial x_N}\right| \le c_{11}\frac{|x_N|}{|x|^N}+c_{12}.
\end{eqnarray*}
Therefore, we have that
$$
\left|\mathbb{G}_{B_1(0)}[  \frac{\partial\delta_{0}}{\partial x_N}]\right|\le \frac{c_{13}}{|x|^{N-1}},\quad \forall x\in B_1(0)\setminus\{0\}.
$$

{\it Proof of (\ref{s53.0}). }  Let $E$ be a Borel set of $B_1(0)$ with $|E|>0$,  then there exists $r_1\in(0,1]$ such that
$$  |E|=|B_{r_1}(0) |.$$
 We deduce  that
\begin{eqnarray*}
  \int_E|\mathbb{G}_{B_1(0)}[  \frac{\partial\delta_{0}}{\partial x_N}]| dx &=& \int_{E\cap B_{r_1}(0)} \frac{c_{13}}{|x|^{N-1}}dx+\int_{E\setminus B_{r_1}(0)} \frac{c_{13}}{|x|^{N-1}}dx\\&\le & \int_{  B_{r_1}(0)} \frac{c_{13}}{|x|^{N-1}}dx
    \\&=& c_{14} r_1= c_{15}\left(\int_E |x|dx\right)^{\frac1{N}}.
\end{eqnarray*}
By the definition of Marcinkiewicz space, we have that
\begin{eqnarray*}
\|\mathbb{G}_{B_1(0)}[  \frac{\partial\delta_{0}}{\partial x_N}]\|_{M^{\frac{N}{N-1}}(B_1(0),\, dx)} \le c_{15}.
\end{eqnarray*}

{\it Proof of (\ref{s53.1}). }  Let $E$ be a Borel set of $B_1(0)$ with $|E|>0$,  then there exists $r_2\in(0,1]$ such that
$$\int_E |x|dx=\int_{B_{r_2}(0)} |x|dx.$$
Since
$$\int_{B_{r_2}(0)} |x|dx=c_{16}r_2^{N+1},$$
 we deduce that
\begin{eqnarray*}
  \int_E|\mathbb{G}_{B_1(0)}[  \frac{\partial\delta_{0}}{\partial x_N}]||x|dx &\le& \int_E \frac{c_{13}}{|x|^{N-1}}|x|dx
    \\&\le & c_{17} r_2^2= c_{18}(\int_E |x|dx)^{\frac2{N+1}}\le c_{19}.
\end{eqnarray*}
 This ends the proof. \qquad$\Box$\smallskip
\begin{lemma}\label{lm 4.1}
Assume that $g:[0,\infty)\to[0,\infty)$ is  continuous,
nondecreasing and verifies (\ref{1.2}). Then for $q\ge \frac{N+1}{N-1}$,
$$\lim_{s\to\infty}g(s)s^{-q}=0.$$

\end{lemma}
{\bf Proof.} Since
\begin{eqnarray*} \int_s^{2s}g(t)t^{-1-k_{\alpha,\beta}}dt\ge
g(s)(2s)^{-1-\frac{N+1}{N-1}}\int_s^{2s}dt=2^{-1-\frac{N+1}{N-1}}g(s)s^{-\frac{N+1}{N-1}}
\end{eqnarray*}
and by (\ref{1.4}),
\begin{eqnarray*}
\lim_{s\to\infty}\int_s^{2s}g(t)t^{-1-\frac{N+1}{N-1}}dt=0,
\end{eqnarray*}
then
$$\lim_{s\to\infty}g(s)s^{-\frac{N+1}{N-1}}=0.$$
The proof is complete.\hfill$\Box$\\[2mm]

Now we are ready to prove Theorem \ref{teo 3.1}.\smallskip

{\bf Proof Of Theorem \ref{teo 3.1}.}
{\it Existence.} Let $\{g_n\}$ be a sequence of $C^{1}$ nondecreasing  functions defined by (\ref{3.1}).
It follows that $\{g_n\}$ is a sequence of odd, bounded and nondecreasing Lipschiz continuous functions.

By Proposition \ref{pr s52.2},  problem (\ref{eq 3.1})
admits a unique $x_N$-odd weak solution $w_{k,n}$ such that
\begin{equation}\label{3.0}
 0\le w_{k,n}\le 2k\mathbb{G}_{B_1(0)}[\frac{\partial  \delta_0}{\partial x_N }]\quad{\rm a.e.\ in}\quad B_1^+(0)
\end{equation}
and
\begin{equation}\label{s52.1.1000}
\int_{B_1(0)} [w_{k,n}(-\Delta) \xi+g_n(w_{k,n})\xi]dx=2k \frac{\partial  \xi(0)}{\partial x_N},\quad \forall\xi\in C_0^{1,1}(B_1(0)).
\end{equation}

For $\beta\in(0,1)$, any compact set $K$ and open set $O$ in $B_1(0)$ satisfying $K\subset O$, $0\not\in\bar O$,
 we have that
$$\norm{w_{k,n}}_{C^{2,\beta}(K)}\le c_{20}\norm{\mathbb{G}_{B_1(0)}[\frac{\partial  \delta_0}{\partial x_N }]}_{C^1(O)},$$
where $c_{20}>0$.
Therefore, up to some subsequence, there exists $w_k$ such that
$$\lim_{n\to0^+}w_{k,n}=w_k\quad{\rm a.e.\ in}\ B_1(0).$$
Then $\{ g_n(w_{k,n})\}$ converges to $g(w_k)$ a.e. in $B_1(0)$.
By Lemma \ref{lm s53.1}, we have that
$$w_{k,n}\to w_k\ {\rm in}\ L^1({B_1(0)})\quad {\rm as}\ \ n\to+\infty,$$
$$\norm{g_n(w_{k,n})}_{L^1(B_1(0),|x|dx)}\le c_{21}\norm{\mathbb{G}_{B_1(0)}[\frac{\partial  \delta_0}{\partial x_N }]}_{L^1(B_1(0),|x|dx)},$$
 by Proposition \ref{pr 1} and $\mathbb{G}_{B_1(0)}[\frac{\partial \delta_0}{\partial x_N}]\in  M^{\frac{N+1}{N-1}}(B_1(0), |x|dx)$, we have that
$$ m(\lambda)\leq c_{22}\lambda^{-\frac{N+1}{N-1}}, \ \quad \forall   \lambda>\lambda_0,$$
where
$$  m(\lambda)=\int_{ S_\lambda}|x|dx
\quad{\rm  with}\quad   S_\lambda=\left\{x\in B_1(0): \left|\mathbb{G}_{B_1(0)}[\frac{\partial \delta_0}{\partial x_N}]\right|>\lambda\right\}.$$
For any Borel
set $E\subset{B_1(0)}$, we have that
\begin{eqnarray*}
 \int_{E}|g_n(w_{k,n})| |x| dx &\le &\int_{E\cap S^c_{\frac{\lambda}{2k}}}g\left(2k\left|\mathbb{G}_{B_1(0)}[\frac{\partial \delta_0}{\partial x_N}]\right|\right)|x| dx
 \\&&+\int_{E\cap S_{\frac{\lambda}{2k}}}g\left(2k\left|\mathbb{G}_{B_1(0)}[\frac{\partial \delta_0}{\partial x_N}]\right|\right)|x|dx \\
   &\le & g(\lambda)\int_E|x|dx+\int_{S_{\frac{\lambda}{2k}}} g\left(2k\left|\mathbb{G}_{B_1(0)}[\frac{\partial \delta_0}{\partial x_N}]\right|\right)|x|dx
   \\&\le&
g(\lambda)\int_E|x|dx+ m\left(\frac{\lambda}{2k}\right)g(\lambda)  +\int_{\frac{\lambda}{2k}}^\infty  m(s)d g(2ks).
\end{eqnarray*}

On the other hand,
$$\int_{\frac{\lambda}{2k}}^\infty   g(2ks)d  m(s)=\lim_{T\to\infty}\int_{\frac{\lambda}{2k}}^{\frac{T}{2k}}  g(2ks)d  m(s).
$$
Thus,
\begin{eqnarray*}
  m\left(\frac{\lambda}{2k}\right)   g (\lambda)+ \int_{\frac{\lambda}{2k}}^{\frac{T}{2k}}   m(s)d  g(2ks) &\le & c_{24}  g(\lambda)\left(\frac{\lambda}{2k}\right)^{-\frac{N+1}{N-1}}+c_{24}\int_{\frac{\lambda}{2k}}^{\frac{T}{2k}} s^{-\frac{N+1}{N-1}}d  g(2ks) \\
   &=&  c_{25}T^{-\frac{N+1}{N-1}}g(T)+c_{26}\int_{\frac{\lambda}{2k}}^{\frac{T}{2k}} s^{-1-\frac{N+1}{N-1}} g(s)ds,
\end{eqnarray*}
where $c_{26}=\frac{c_{24}}{\frac{N+1}{N-1}+1}(2k)^{2+\frac{N+1}{N-1}}$.
By assumption (\ref{1.2}) and Lemma \ref{lm 4.1}, we have that  $T^{-\frac{N+1}{N-1}}  g(T)\to 0$ as $T\to\infty$, therefore,
$$ m\left(\frac{\lambda}{2k}\right)   g ( \lambda )+ \int_{\frac{\lambda}{2k}}^\infty   m(s)\ d  g(2ks)\leq c_{26}\int_{\frac{\lambda}{k}}^\infty s^{-1-\frac{N+1}{N-1}}  g(s)ds.
$$
Notice that the  quantity on the right-hand side tends to $0$
when $\lambda\to\infty$. The conclusion follows: for any
$\epsilon>0$, there exists $\lambda>0$ such that
$$c_{26}\int_{\frac{\lambda}{2k}}^\infty s^{-1-\frac{N+1}{N-1}}  g(s)ds\leq \frac{\epsilon}{2}.
$$
For $\lambda$ fixed,  there exists $\delta>0$ such that
$$\int_E|x| dx\leq \delta\Longrightarrow  g ( \lambda)\int_E|x| dx\leq\frac{\epsilon}{2},
$$
which implies that $\{g_n\circ w_{k,n}\}$ is uniformly integrable in
$L^1(B_1(0),|x| dx)$. Then $g_n\circ w_{k,n}\to g\circ w_k$ in
$L^1(B_1(0),|x| dx)$ by Vitali convergence theorem, see \cite{G}.

Furthermore, for any $\xi\in C^{1,1}_0(B_1(0))$, we know that
\begin{equation}\label{3.01}
 |\xi(x)-\xi(0)-\nabla \xi(0)\cdot x|\le c_{27}|x|^2,
\end{equation}
and it follows the odd prosperity of $w_{k,n}$, $g_n$ and $g$,
  that
\begin{equation}\label{3.02}
 \int_{B_1(0)}g_n(w_{k,n})\xi(0)\,dx=0\quad{\rm and}\quad \int_{B_1(0)}g(w_{k})\xi(0)\,dx=0,
\end{equation}
then
\begin{eqnarray*}
&&\left|\int_{B_1(0)}g_n(w_{k,n})\xi\,dx-\int_{B_1(0) }g(w_{k})\xi\,dx\right|
\\&\le&\left|\int_{B_1(0)}[g_n(w_{k,n})-g(w_{k})]\nabla \xi(0)\cdot x\,dx\right| +c_{27} \int_{B_1(0)}|g_n(w_{k,n})-g(w_{k})| |x|^2\,dx   \\
   &\le &  c_{28}\norm{g_n(w_{k,n})-g(w_{k})}_{L^1(B_1(0),\,|x| dx)}\to0\quad {\rm as}\quad n\to+\infty.
\end{eqnarray*}

Then passing to the limit as
$n\to +\infty$ in the identity (\ref{s52.1.1000}),
it implies that
\begin{equation}\label{weak 1}
 \int_{B_1(0)} [w_k(-\Delta)\xi+g(w_k)\xi]dx=2k\frac{\partial \xi(0)}{\partial x_N}.
\end{equation}
Thus, $w_{k}$ is a very weak solution of (\ref{eq 1.1}). The regularity results follows by Lemma \ref{lm 3.1}.

{\it Proof of $(i)$.} Since $w_{k,n}$ is $x_N$-odd in $B_1(0)\setminus\{0\}$ and $w_k=\lim_{n\to+\infty} w_{k,n}$ in $B_1(0)\setminus\{0\}$,
then it implies that $w_k$ is $x_N$-odd in $B_1(0)\setminus\{0\}$. By the fact that $w_{k+1,n}\ge w_{k,n}\ge0 $ in $B_1^+(0)$, it follows  that
$$w_{k+1} \ge w_{k}\ge0\quad {\rm in}\ B_1^+(0).$$

{\it Proof of $(ii)$.}
We observe that
$$0\le w_{k,n}\le  2k\mathbb{G}_{B_1(0)}[\frac{\partial \delta_0}{\partial x_N}]\quad {\rm in}\ B_1^+(0),$$
then $g_n(w_{k,n})\le g(2k\mathbb{G}_{B_1(0)}[\frac{\partial \delta_0}{\partial x_N}])$ for $x\in B_1^+(0)$ and
let $w_g$ be the unique solution of (\ref{homo}) with $f=g(2k\mathbb{G}_{B_1(0)}[\frac{\partial \delta_0}{\partial x_N}])$, we have that
  \begin{equation}\label{s53.3}
 2k\mathbb{G}_{B_1(0)}[\frac{\partial \delta_0}{\partial x_N}](x) \ge  w_{k,n}(x)  \ge  2k\mathbb{G}_{B_1(0)}[\frac{\partial \delta_0}{\partial x_N}](x)-w_g(x).
 \end{equation}
Then $w_k$ satisfies (\ref{s53.3}).
From Proposition \ref{pr 2.1}, it infers that for $x=te$ with $t\in(0,\frac12)$ $e=(e_1,\cdots,e_N)\in\partial B_1(0)$ with $  e_N>0$,
\begin{eqnarray*}
w_g(te)&  \le & \int_{B_1^+(0)} \left[\frac{|y|}{|y-te| (|y-te|+2|y|)^{N-2}}\right]g\left( 2k|\mathbb{G}_{B_1(0)}[\frac{\partial \delta_0}{\partial x_N}](y)|\right)dy \\
&\le &\int_{B_1^+(0)} \left[\frac{|y|}{|y-te| (|y-te|+2|y|)^{N-2}}\right]g(2c_{13}k|y|^{1-N})dy\\
 &    := & \int_{B_1^+(0)} A(t,y) dy.
\end{eqnarray*}

For $y\in B_{\frac t2}(te)$, we have that $|y-te|\ge \frac t2$ and $\frac t2\le |y|\le \frac{3t}2$,
then
\begin{eqnarray*}
t^{N-1}\int_{B_{\frac t2}(te)}A(t,y) dy\le c_{30} r^{-\frac{N+1}{N-1}}g ( 2c_{13}k r )
 \to0\quad{{\rm as}}\ \  r\to+\infty,
\end{eqnarray*}
where $r=t^{-\frac1{N-1}}$.

For $y\in B_{\frac t2}(0)$, we have that $|y-te|\ge \frac t2$
and
\begin{eqnarray*}
t^{N-1}\int_{B^+_{\frac t2}(0)}A(t,y) dy&\le & \int_{B_{\frac t2}(0)}|y|g(2c_{13}k|y|^{1-N})dy
\\&=& c_{32}\int_0^{\frac t2}g(2c_{13}ks^{1-N})s^{N}ds
\\ &=&c_{32}\int_{(\frac t2)^{-\frac1{N-1}}}^{\infty}  g ( 2c_{13}k \tau )\tau^{-1-\frac{N+1}{N-1} }d\tau
\\&\to&0\quad{{\rm as}}\ \ t\to0,
\end{eqnarray*}
where we have used (\ref{1.2}).

For $y\in B_1^+(0)\setminus (B_{\frac t2}(0)\cup B_{\frac t2}(te))$, we have that
$|y-te|\ge \frac t2, |y|\ge \frac t2$
and
\begin{eqnarray*}
t^{N-1}\int_{B_1^+(0)\setminus (B_{\frac t2}(0)\cup B_{\frac t2}(te))}A(t,y) dy&\le & c_{33}t\int_{B_1^+(0)\setminus (B_{\frac t2}(0)\cup B_{\frac t2}(te))}  g(2c_{13}k|y|^{1-N})dy
\\&=& c_{33}t \int_{\frac t2}^1g(2c_{13}kr^{1-N}) r^{N-1}dr
\\ &=&c_{34}\frac{t^{N-1}g(2c_{13}kt^{1-N})}{t^{-2}}=c_{35}\tau^{-\frac{N+1}{N-1}}g(\tau)
\\&\to&0\quad{{\rm as}}\ \ \tau\to+\infty
\end{eqnarray*}
where $\tau=2c_{13}kt^{1-N}$.
 Then  for any $e=(e_1,\cdots,e_N)\in \partial B_1(0)$ and $e_N\not=0$, we have that
 \begin{equation}\label{3.001}
 \lim_{t\to0^+}t^{N-1}w_g(te)=0
 \end{equation}
 and
$$
\lim_{t\to0^+}\mathbb{G}_{B_1(0)}[g(\mathbb{G}_{B_1(0)}[2k\frac{\partial  \delta_{0}}{\partial x_N}])](te)t^{N-1}=0.
$$
 By (\ref{s53.3}),
$$\lim_{t\to0^+} \mathbb{G}_{B_1(0)}[\frac{\partial  \delta_{0}}{\partial x_N}](te)t^{N-1}=  e_N, $$
 and $x_N$-odd property of $w_{k}$,
we derive that
 for any $e\in \partial B_1(0)$,
$$\lim_{t\to0^+} w_k(te)t^{N-1}=2k e_N. $$

Finally, we prove the uniqueness. Let $u_k,\, v_k$ be two $x_N-$odd solutions of (\ref{eq 1.1}), from Lemma \ref{lm 3.3},  $u_k,\, v_k$ are two solutions of
\begin{equation}\label{2.02}
\arraycolsep=1pt\left\{
\begin{array}{lll}
 -\Delta    u+g(u)=0\quad  &{\rm in}\quad\ \ B_1^+(0),\\[2mm]
 \phantom{-\Delta +g(u) }
 u=k\delta_0\quad &{\rm on}\quad\ \  \partial B_1^+(0).
\end{array}\right.
\end{equation}
Form  the uniqueness of the very weak solution to (\ref{2.02}), see  \cite[Theorem 2.1]{GV}, we obtain $u_k=v_k$ in $B_1^+(0)$ and combine the  $x_N-$odd property we obtain
the uniqueness of the very weak solution of (\ref{eq 1.1}) with $j=0$.
This ends the proof.\qquad$\Box$

\subsection{  Nonexistence  }

This subsection is devoted to  obtain the nonexistence of very weak solutions of (\ref{eq 1.1}) in the supercritical case.

\begin{lemma}\label{lm 3.3}
Assume that $k>0$,   the nonlinearity $g:\R\to\R$ is an odd, nondecreasing and Lipchitz continuous function.
Let $u$ be an $x_N$-odd very weak solution of  (\ref{eq 1.1}).

 Then $u$ is a very weak solution of (\ref{2.02}).

\end{lemma}
{\bf Proof.} From Lemma \ref{lm 3.1} and the $x_N$-odd property,
we have that
$$u=0\quad{\rm on}\quad \partial B_1^+(0)\setminus\{0\}.$$
For any $x_N$-odd function $\xi\in C^{1,1}_0(B_1(0))$,
we deduce that
$$\xi=0\quad{\rm on}\quad \partial B_1^+(0)\qquad |\xi(x)|\le c_{36}|x|,\quad x\in B_1(0),$$
where  $c_{36}>0$.
Since $g(u)\in L^1(B_1(0),|x|dx)$,  then for $x_N$-odd function $\xi\in C^{1,1}_0(B_1(0))$, we have that
$g(u)\xi\in L^1(B_1(0))$ and the weak solution $u$ satisfies that
\begin{equation}\label{2.01}
  \int_{B_1(0)}[ u(-\Delta)\xi dx+   g(u)\xi ]dx=2k\frac{\partial \xi(0)}{\partial x_N},
\end{equation}
for   $x_N-$odd function $\xi\ {\rm in}\ C^{1,1}_{0}(B_1(0))$.
By $x_N$-odd property, we have that
$$\int_{B_1^+(0)}[ u(-\Delta)\xi dx+   g(u)\xi ]dx=\int_{B_1^-(0)}[ u(-\Delta)\xi dx+   g(u)\xi ]dx,$$
which, combined with (\ref{2.01}), implies that
$$
  \int_{B_1^+(0)}[ u(-\Delta)\xi dx+   g(u)\xi ]dx=k\frac{\partial \xi(0)}{\partial x_N},\qquad \forall \xi\in C^{1,1}_{0}(B_1^+(0)).
$$
So we have that $u$ is a weak solution of (\ref{2.02}).\qquad $\Box$
\medskip

\noindent{\bf Proof of Theorem \ref{teo unique}.}
When $g(s)=|s|^{p-1}s$ with $p\ge \frac{N+1}{N-1}$ and $k>0$, \cite[Theorem 3.1]{GV} shows that the nonnegative solution of
$$
\arraycolsep=1pt\left\{
\begin{array}{lll}
 -\Delta    u+g(u)=0\quad  &{\rm in}\quad\ \ B_1^+(0),\\[2mm]
 \phantom{-\Delta +g(u) }
 u=0\quad &{\rm on}\quad\ \  \partial B_1^+(0)\setminus\{0\}
\end{array}\right.
$$
has removable singularity at origin, so  problem (\ref{2.02}) has no very weak solution, which contradicts Lemma \ref{lm 3.3}.
Then (\ref{eq 1.1}) has no $x_N$-odd weak solution when $g(s)=|s|^{p-1}s$ with $p\ge \frac{N+1}{N-1}$.\qquad$\Box$

\setcounter{equation}{0}

\section{Strongly anisotropic singularity for $p\in(1,\frac{N+1}{N-1})$}

In this section, we consider the limit of weak solutions $w_k$ as $k\to\infty$ to
\begin{equation}\label{s54.3}
\arraycolsep=1pt\left\{
\begin{array}{lll}
  -\Delta    u+|u|^{p-1}u=2k\frac{\partial \delta_0}{\partial x_N}\quad  &{\rm in}\quad\ \ B_1(0),\\[2mm]
 \phantom{-\Delta    u+|u|^{p-1} }
 u=0\quad &{\rm on}\quad\ \ \partial B_1(0),
 \end{array}\right.
\end{equation}
 where   $p\in(1,\frac{N+1}{N-1})$.
By Theorem \ref{teo 1}, we observe  that the mapping $k\mapsto w_k$ is   nondecreasing in $B_1^+(0)$,
then $\lim_{k\to+\infty} w_k(x)$ exists for any $x\in B_1(0)\setminus\{0\}$,  denoting
\begin{equation}\label{4.1}
 w_\infty(x)=\lim_{k\to+\infty} w_k(x)\quad{\rm for}\ \ x\in B_1(0)\setminus\{0\}.
\end{equation}

For $w_\infty$, we have the following result.

\begin{proposition}\label{pr s54.1}
Let $p\in(1,\frac{N+1}{N-1})$, then $w_\infty$ is $x_N$-odd,
$$0\le w_\infty(x)\le \lambda_0 |x|^{-\frac2{p-1}},\quad\forall\, x\in  B_1^+(0)$$
for some $\lambda_0>0$
and $w_\infty$ is a classical solution of
\begin{equation}\label{s54.1}
\arraycolsep=1pt\left\{
\begin{array}{lll}
 -\Delta    u+|u|^{p-1}u=0\quad  &{\rm in}\quad \ B_1(0)\setminus\{0\},\\[2mm]
 \phantom{-\Delta    u+|u|^{p-1} }
 u=0\quad &{\rm on}\quad\ \ \partial B_1(0).
 \end{array}\right.
\end{equation}

\end{proposition}
{\bf Proof.} In order to  obtain the asymptotic behavior of $w_\infty$ near the origin,
 we construct the function
$$
  v_p(x)=|x|^{-\frac{2}{p-1}},\qquad\forall\, x\in\R^N\setminus\{0\}.
$$
For $p\in(1,\frac{N+1}{N-1})$, there exists $\lambda_0>0$ such that $\lambda_0 v_p$ is a super solution of
(\ref{s54.1}) and $-\lambda_0 v_p$ is a sub solution of (\ref{s54.1}).

It follows by Theorem \ref{teo 1} that  $w_k$ is a classical solution of (\ref{s54.1}) satisfying (\ref{1.3}), hence by Comparison Principle, for any $k$, there exists $r_k\in(0,1)$ small enough such that
$$|w_k|\le \lambda_0 v_p\quad {\rm in}\quad   B_{r_k}(0)\setminus \{0\}.$$
By Comparison Principle, we have that
$$|w_k|\le \lambda_0 v_p\quad {\rm in}\quad B_1(0) \setminus \{0\}.$$
Since $k$ is  arbitrary, we deduce that
$$|w_\infty|\le \lambda_0 v_p\quad {\rm in}\quad B_1(0) \setminus \{0\}.$$
Therefore, from standard Stability Theorem, we derive that $w_\infty$ is a classical solution of  (\ref{s54.1}).
\qquad $\Box$

\smallskip

We next do a precise bound for $w_\infty$ to prove (\ref{1.5}).
\begin{lemma}\label{lm 5.1}
Let $p\in(1,\frac{N+1}{N-1})$ and   $w_\infty$ be defined by (\ref{4.1}), then $w_\infty$ satisfies (\ref{1.5}).

\end{lemma}
{\bf Proof.}
{\it We claim that
\begin{equation}\label{4.3}
 \frac1{c_{36}}t^{\frac2{p-1}}P_N(e)\le w_\infty(te)\le c_{36}t^{\frac2{p-1}}P_N(e),\quad \forall t>0,\ e=(e_1,\cdots,e_N)\in \partial B_1(0),  e_N>0.
\end{equation}
}
We have that
\begin{equation}\label{4.2}
 2k\mathbb{G}_{B_1(0)}[\frac{\partial \delta_0}{\partial x_N}](x) \ge  w_{k}(x)  \ge  2k\mathbb{G}_{B_1(0)}[\frac{\partial \delta_0}{\partial x_N}](x)-\varphi_p(x),\quad\forall x\in B^+_1(0),
 \end{equation}
where $\varphi_p:=k^p\mathbb{G}_{\R^N}[\mathbb{G}_{\R^N}[\frac{\partial \delta_0}{\partial x_N}]^p]  \ge  \mathbb{G}_{B_1(0)}[w_k^p]$, and $\varphi_p$ is a $x_N$-odd  solution of
\begin{equation}\label{eq 4.1}
\arraycolsep=1pt\left\{
\begin{array}{lll}
 \displaystyle  \ -\Delta    u =\tilde c_N^p \frac{ |x_N|^{p-1}x_N }{|x|^{Np} }\quad  &{\rm in}\quad\ \ \R^N\setminus \{0\},\\[2.5mm]
 \phantom{ }
 \displaystyle  \lim_{|x|\to+\infty}u(x)=0.
 \end{array}\right.
\end{equation}
Indeed, $\mathbb{G}_{\R^N}[\mathbb{G}_{\R^N}[\frac{\partial \delta_0}{\partial x_N}]^p]$ is $x_N-$odd and
\begin{eqnarray*}
 \varphi_p(x) &=& c_N\tilde c_N^p\int_{\R^N}\frac{|y_N|^{p-1}y_N}{|x-y|^{N-2}|y|^{Np}} dy  \\
   &=&c_N\tilde c_N^p|x|^{(1-N)p+2} \int_{\R^N}\frac{|y_N|^{p-1}y_N}{|e_z-y|^{N-2}|y|^{Np}} dx,
\end{eqnarray*}
where $e_z=\frac{z}{|z|}$.

We observe that    $\varphi_p$ satisfies
$$ -\Delta    u(x) =\tilde c_N^p \frac{ |x_N|^{p-1}x_N }{|x|^{Np} },\quad     \forall  x\in B_{\frac32}^+(0)\setminus B_{\frac12}^+(0), $$
and then it follows by Hopf's Lemma (see \cite{E}) that
$$\varphi_p(x)\le c_{37} x_N,\quad\forall x\in \partial B_1^+(0).$$
Therefore,
$$\varphi_p(e)\le c_{38}P_N(e),\quad\forall e\in \partial B_1^+(0).$$

  {\it Proof of lower bound in (\ref{4.3}).} It follows by (\ref{4.2})  that
$$w_k(x)\ge c_{40} k|x|^{1-N}P_N(\frac x{|x|})-c_{39}k^p|x|^{(1-N)p+2}P_N(\frac x{|x|}),\quad\forall\, x\in B_{\frac12}(0)\setminus\{0\}.$$
Set $\rho_k= (2^{(N-1)(p-1)-3}\frac{c_{39}}{c_{40}}k^{p-1} )^{\frac1{(N-1)(p-1)-2}}$,  then
\begin{eqnarray*}
 c_{39}k^p|x|^{(1-N)p+2} &\le &  c_{39}k^p(\frac{\rho_k}2)^{(1-N)p+2} \\
   &\le &  \frac{c_{40}} k \rho_{k}^{1-N}\le  \frac{c_{40}}2 k |x|^{1-N}
\end{eqnarray*}
and
$$k=c_{41} \rho_k ^{1-N-\frac2{p-1}}\ge c_{41}|x|^{N-1-\frac2{p-1}},$$
where $c_{41}>0$ independent of $k$.
Thus, for $t\in(\frac{\rho_k}2,\rho_k),\ e\in \partial B_1(0),\ e\cdot e_N>0,$
\begin{eqnarray*}
 w_k(te) &\ge & [c_{40}kt^{1-N}-c_{39}k^pt^{(1-N)p+2}]P_N(e)\\
   &\ge &   \frac{c_{40}}2 k t^{1-N}P_N(e)\\
    &\ge & \frac{c_{40}c_{41}}2 t^{-\frac2{p-1}}P_N(e).
\end{eqnarray*}
Now we can choose a sequence $\{k_n\}\subset [1,\infty)$ such that
$$\rho_{k_{n+1}}\ge \frac12 \rho_{k_n}$$
and for any $x\in B_{\frac12}^+(0)\setminus\{0\}$, there exists $k_n$ such that $x\in B_{\rho_k}^+(0)\setminus B_{\frac{\rho_k}2}^+(0)$
and
$$w_{k_n}(x)\ge \frac{c_{40}c_{41}}2|x|^{-\frac2{p-1}}P_N(\frac{x}{|x|}).$$
Together with $w_{k_{n+1}}\ge w_{k_n}$ in $B_1^+(0)$, we have that
$$w_\infty(x)\ge \frac{c_{40}c_{41}}2|x|^{-\frac2{p-1}}P_N(\frac{x}{|x|}),\quad\forall\, x\in B_1^+(0).$$

{\it Proof of the upper bound in (\ref{4.3}).} Let $\bar w_p=|x|^{-\frac2{p-1}-1}x_N$, then
$$-\Delta \bar w_p=(\frac2{p-1}+1)(\frac2{p-1}+1-N)|x|^{-\frac2{p-1}-3}x_N, $$
where $(\frac2{p-1}+1)(\frac2{p-1}+1-N)>0$ for $p<\frac{N+1}{N-1}$.
By Comparison Principle, there exists $t_0>0$ independent of $k$ such that
$$u_k\le t_0 \bar w_p\quad{\rm in}\quad B_1^+(0), $$
which implies that
$$w_\infty \le t_0\bar w_p\quad{\rm in}\ \ B_1^+(0).$$

{\it Proof of (\ref{1.5}).}  We observe  there exists $t_{00}\in(0,t_0)$ such that
$$-\Delta t_{00}\bar w_p+(t_{00}\bar w_p)^p\le 0\quad{\rm in}\ \ \R^N_+.$$
Therefore,
$$\arraycolsep=1pt\left\{
\begin{array}{lll}
 -\Delta u_p+u_p^p= 0\quad  &{\rm in}\quad\ \ \R^N_+,\\[2mm]
 \phantom{----\  }
  u=0\quad  &{\rm on}\quad\ \ \partial\R^N_+\setminus\{0\}
 \end{array}\right.
 $$
admits a unique solution $u_p$ and by scaling property, we have that
$$u_p(te)=t^{-\frac2{p-1}}u_p(e),\quad te\in \R^N_+.$$
By Comparison Principle, we have that
$$ u_p(te)-\max u_p(e)\le  u_\infty(te)\le u_p(te),\quad te\in B_1^+(0),$$
which implies (\ref{1.5}) with $\varphi(e)=u_p(e)$.
\qquad$\Box$

\medskip

\noindent{\bf Proof of Theorem \ref{teo 3}.} From Proposition \ref{pr s54.1} and Lemma \ref{lm 5.1},
one has that $u_\infty$ is a classical solution of (\ref{eq 1.03}) satisfying (\ref{1.5}).\qquad$\Box$

\setcounter{equation}{0}

\section{Non $x_N$-odd solutions }

\subsection{Existence}

 Under the assumptions on $g$ in Theorem 1.1, it shows  from \cite{V} that the problem
\begin{equation}\label{eq 5.1}
\arraycolsep=1pt\left\{
\begin{array}{lll}
 -\Delta    u+g(u)=j\delta_0\quad  &{\rm in}\quad\ \ B_1(0),\\[2mm]
 \phantom{-\Delta +g(u) }
 u=0\quad &{\rm on}\quad\ \  \partial B_1(0)
\end{array}\right.
\end{equation}
admits  a unique weak solution, denoting by $u_{0,j}$. In the approaching the weak solution of problem (\ref{eq 1.1}) with $j>0$,
 a barrier will be constructed by $u_{0,j}$ and $u_{k}$, where $u_{k}$ is the unique $x_N-$odd weak solution of (\ref{eq 1.1}) with $j=0$.

\medskip

\noindent{\bf Proof of Theorem \ref{teo 1} with $j> 0$. }
{\it Step 1.} We observe that $g_n$ is bounded, Lipschitz continuous and nondecreasing, where   $g_n$ is defined by (\ref{3.1}),
 then it follows from \cite[Theorem  3.7]{V} and the Kato's inequality  that
 \begin{equation}\label{eq 5.2}
\arraycolsep=1pt\left\{
\begin{array}{lll}
 -\Delta    u+g_n(u)=k\mu_t+j\delta_0\quad  &{\rm in}\quad\ \ B_1(0),\\[2mm]
 \phantom{-\Delta +g_n(u) }
 u=0\quad &{\rm on}\quad\ \  \partial B_1(0)
 \end{array}\right.
\end{equation}
  admits a unique weak solution $v_{k,j,n,t}$, which is a classical solution of
$$
\arraycolsep=1pt\left\{
\begin{array}{lll}
 -\Delta    u+g_n(u)=0\quad  &{\rm in}\quad\ \ B_1(0)\setminus\{te_N,0,-te_N\},\\[2mm]
 \phantom{-\Delta +g_n(u) }
 u=0\quad &{\rm on}\quad\ \  \partial B_1(0).
 \end{array}\right.
$$ Moreover, $v_{k,j,n,t}$ could be approximated by
the classical solutions $\{v_{n,t,m}\}$ to
$$
\arraycolsep=1pt\left\{
\begin{array}{lll}
 -\Delta    u+g_n(u)=k\mu_{t,m}+j\sigma_m\quad  &{\rm in}\quad\ \ B_1(0),\\[2mm]
 \phantom{-\Delta +g_n(u) }
 u=0\quad &{\rm on}\quad\ \  \partial B_1(0),
 \end{array}\right.
 $$
where $$\mu_{t,m}(x)=\frac{\sigma_m(x-te_N)-\sigma_m(x+te_N)}t$$ and
$\{\sigma_m\}$ is a sequence of radially symmetric, nondecreasing smooth functions converging
to $\delta_0$ in the distribution sense. Furthermore,
\begin{equation}\label{5.3}
\int_{B_1(0)} [v_{n,t,m}(-\Delta)\xi+ g_n(v_{n,t,m})\xi]\, dx=\int_{B_1(0)}[k\mu_{t,m}+j\sigma_m]\xi dx,\quad \forall \xi\in C^{1,1}_0(B_1(0)).
\end{equation}

Since  $k\mu_{t,m}+j\sigma_m \ge k\mu_{t,m}$, it implies by  Comparison Prinsiple that
\begin{equation}\label{5.5}
w_{k,n,t,m}\le v_{n,t,m}\le w_{k,n,t,m}+\upsilon_{j,m}\quad {\rm in}\ B_1(0),
\end{equation}
where $w_{k,n,t,m}$ is the weak solution of (\ref{3.2}) and
$\upsilon_{j,n,m}$ is the unique solution of the   equation
$$
\arraycolsep=1pt\left\{
\begin{array}{lll}
 -\Delta    u =j\sigma_m\quad  &{\rm in}\quad\ \ B_1(0),\\[2mm]
 \phantom{-\Delta  }
 u=0\quad &{\rm on}\quad\ \  \partial B_1(0).
 \end{array}\right.
 $$
Therefore,   $v_{k,j,n,t}$ satisfies
 $$w_{k,n,t}\le v_{k,j,n,t}\le w_{k,n,t}+ j\upsilon_{j,m}\quad {\rm in}\quad\ \ B_1(0)\setminus\{te_N,0,-te_N\}$$
and
\begin{equation}\label{5.2}
 v_{k',j',n,t}\ge v_{k,j,n,t}  \quad {\rm in}\quad B_1^+(0)\setminus\{te_N\}\quad {\rm for}\ \ k'\ge k,\ j'\ge j.
\end{equation}

{\it Step 2.} From Step 1,  problem (\ref{eq 5.2}) admits a unique weak solution $v_{k,j,n,t}$,
that is,
\begin{equation}\label{5.4}
\int_{B_1(0)} [w_{k,n,t}(-\Delta)\xi+g_n(w_{k,n,t})\xi ]\, dx=\frac{\xi(te_N)-\xi(-te_N)}{t}+j\xi(0),\quad \forall \xi\in  C^2_0(B_1(0)).
\end{equation}
On the one hand, by  Lemma \ref{lm 2.1}, we have that
$$\lim_{t\to0^+}\frac{\xi(te_N)-\xi(-te_N)}{t}=2\frac{\partial\xi(0)}{\partial x_N}.$$
On the other hand, by the fact that
$$|w_{k,n,t}|\le k|\mathbb{G}_{B_1(0)}[\mu_t]|\quad{\rm in}\ B_1(0),$$
we have that
$$|v_{k,j,n,t}|\le k|\mathbb{G}_{B_1(0)}[\mu_t]|+j|\mathbb{G}_{B_1(0)}[\delta_0]|\quad{\rm in}\ B_1(0),$$
By interior regularity results, see \cite{HL}, for $\sigma\in(0,1)$ and any compact set $K$ and open set $O$ in $B_1(0)$ such that $K\subset O$, $\bar O\cap \{te_N:\ t\in(-\frac12,\frac12)\}=\emptyset$,
there exist $c_{43},c_{44}>0$ independent of $t$ such that
\begin{eqnarray*}
 \norm{v_{k,j,n,t}}_{C^{2+\sigma}(K)} &\le & c_{43}[\norm{k\mathbb{G}_{B_1(0)}[\mu_{t}] }_{L^\infty(O)}+\norm{j\mathbb{G}_{B_1(0)}[\delta_0] }_{L^\infty(O)}]
  \\ &\le&  c_{44}[k\norm{\mathbb{G}_{B_1(0)}[\frac{\partial\delta_{0}}{\partial x_N}]}_{L^\infty(O)}+\norm{j\mathbb{G}_{B_1(0)}[\delta_0] }_{L^\infty(O)}].
\end{eqnarray*}
Moreover, by  \cite[Lemma 3.6]{CW} $\{\mathbb{G}_{B_1(0)}[\mu_t]\}$ is uniformly bounded in $M^{\frac{N}{N-1}}(B_1(0),dx)$ if $N\ge 3$ and in $M^{\frac{2}{1+\sigma}}(B_1(0),dx)$ for any $\sigma\in(0,\frac12)$ if $N=2$,
$\{\mathbb{G}_{B_1(0)}[\delta_0]\}$ is uniformly bounded in $M^{\frac{N}{N-2}}(B_1(0))$ if $N\ge 3$ and $\{\mathbb{G}_{B_1(0)}[\delta_0]\}$,
 is uniformly bounded in  $M^{q}(B_1(0),dx)$ for any $q>0$.
 Therefore, $\{v_{k,j,n,t} \}_t$ is relatively compact in $L^p({B_1(0)})$ for any $p\in[1,\frac{N}{N-1})$.
Then there exists $v_{k,j,n}\in L^1(B_1(0))$ such that
$$v_{k,j,n,t}\to v_{k,j,n}\quad {\rm  a.e.\ in}\ B_1(0)\quad {\rm and\  in}\  L^1(B_1(0)),$$
which implies that
$$g_n(v_{k,j,n,t})\to g_n(v_{k,j,n})\quad {\rm  a.e.\ in}\ B_1(0)\quad {\rm and\  in}\  L^1(B_1(0)).$$

 Therefore,  up
to some subsequence, passing to the limit as $t\to0^+$ in the identity (\ref{5.3}), it infers that
$v_{k,j,n}$ is the unique very weak solution of
$$
\arraycolsep=1pt\left\{
\begin{array}{lll}
 -\Delta    u+g_n(u)=2k\frac{\partial \delta_0}{\partial x_N}+j\delta_0\quad  &{\rm in}\quad\ \ B_1(0),\\[2mm]
 \phantom{-\Delta +g_n(u) }
 u=0\quad &{\rm on}\quad\ \  \partial B_1(0).
 \end{array}\right.
$$
Here the uniqueness follows by the Kato's inequality.  It follows by (\ref{5.5}), (\ref{3.0}) and $x_N$-odd property of $w_{k,n}$ that
\begin{equation}\label{5.6}
w_{k,n}\le v_{k,j,n}\le w_{k,n}+j\mathbb{G}_{B_1(0)}[\delta_0]\quad {\rm in}\ B_1(0),
\end{equation}
where $w_{k,n}$ is the unique $x_N$-odd solution of (\ref{eq 3.1}). From the proof of Theorem \ref{teo 3.1}, we known that for any $\xi\in C^{1.1}_0(B_1(0))$,
$$\int_{B_1(0)} g_n(w_{k,n})\xi \,dx\to\int_{B_1(0)} g(w_{k})\xi \,dx\quad{\rm as}\quad n\to+\infty. $$

{\it Step 3. }
It follows by (\ref{5.6}) that
\begin{equation}\label{5.7}
\int_{B_1(0)} [v_{k,j,n}(-\Delta) \xi+g_n(v_{k,j,n})\xi]dx=2k \frac{\partial  \xi(0)}{\partial x_N}+j\xi(0),\quad \forall\xi\in C_0^{1,1}(B_1(0)).
\end{equation}

For $\beta\in(0,1)$, any compact set $K$ and open set $O$ in $B_1(0)$ satisfying $K\subset O$, $0\not\in\bar O$,
 we have that
$$\norm{v_{k,j,n}}_{C^{2,\beta}(K)}\le c_{45}[\norm{\mathbb{G}_{B_1(0)}[\frac{\partial  \delta_0}{\partial x_N }]}_{C^1(O)}+\norm{\mathbb{G}_{B_1(0)}[  \delta_0 ]}_{C^1(O)}].$$
Therefore, up to some subsequence, there exists $v_{k,j}$ such that
$$\lim_{n\to+\infty}v_{k,j,n}=v_{k,j}\quad{\rm a.e.\ in}\ B_1(0).$$
Then $\{ g_n(v_{k,j,n})\}$ converges to $g(v_{k,j})$ a.e. in $B_1(0)$.

We observe that  $\tilde v_n:=v_{k,j,n}-w_{k,n}$ is  the very weak solution of
\begin{equation}\label{6.2}
  \arraycolsep=1pt\left\{
\begin{array}{lll}
 -\Delta    u+g_n(w_{k,n}+u)-g_n(w_{k,n})=j\delta_0\quad  &{\rm in}\quad\ \ B_1(0),\\[2mm]
 \phantom{-\Delta +g_n(w_{k,n}+u)-g_n(w_{k,n})}
 u=0\quad &{\rm on}\quad\ \  \partial B_1(0).
 \end{array}\right.
\end{equation}
Note that $0\le \tilde v_n\le  j\mathbb{G}_{B_1(0)}[\delta_0]$ and by (\ref{1.20}), it follows that
$$ 0\le g_n(w_{k,n}+\tilde v_n)-g_n(w_{k,n})\le c_1 \left[\frac{g(w_{k,n})}{1+|w_{k,n}|}\tilde v_n+ g_n(\tilde v_n)\right].$$
Thus, $\{ g_n(w_{k,n}+\tilde v_n)-g_n(w_{k,n})\}$ converges to $g(w_{k}+\tilde v)-g(w_{k})$ a.e. in $B_1(0)$,
where $\tilde v=v_{k,j}-w_k$.

Let $\tilde g_n(s)=g_n(w_{k,n}+s)-g_n(w_{k})$, we see that $\tilde g_n(0)=0$ and function $\tilde g_n$ is nondecreasing and verifies (\ref{1.2}).
Then it follows by Theorem 3.7 in \cite{V}  that
 $$\norm{\tilde g_n  (\tilde v_n)}_{L^1(B_1(0))}\le c_{46}\norm{\mathbb{G}_{B_1(0)}[ \delta_0 ]}_{L^1(B_1(0) dx)}$$
and
$$\tilde v_n\le j\mathbb{G}_{B_1(0)}[\delta_0].$$
Thus, it follows by (\ref{1.20}) that
\begin{eqnarray*}
 \tilde g_n(\tilde v_n(x))  &\le &c_1 g(j\mathbb{G}_{B_1(0)}[ \delta_0 ](x))+c_1\frac{g(2k|\mathbb{G}_{B_1(0)}[ \frac{\partial \delta_0}{\partial x_N} ]|(x))}{1+2k|\mathbb{G}_{B_1(0)}[ \frac{\partial \delta_0}{\partial x_N} ](x)|}j\mathbb{G}_{B_1(0)}[ \delta_0 ](x) \\
   &\le & c_1 g(c_{45}|x|^{2-N}) +c_{46} g(c_{45}|x|^{1-N})|x|.
\end{eqnarray*}

Let $ S_\lambda=\{x\in B_1(0): |x|> \lambda^{-1}\}$, then
$$\tilde  m(\lambda)=\int_{ S_\lambda}dx= c_{47}\lambda^{-N}, \ \quad \forall   \lambda>\lambda_0.$$
For any Borel
set $E\subset{B_1(0)}$, we have that
\begin{eqnarray*}
 \int_{E}|\tilde g_n  (\tilde v_n)|   dx &\le &c_1\int_{E\cap S^c_{ \lambda }}g (c_{45}\lambda^{N-2}) dx +c_1\int_{E\cap S_\lambda }g(c_{45}|x|^{2-N})  dx
 \\&&+c_{46}\int_{E\cap S^c_{ \lambda }}g (c_{45}\lambda^{N-1})\lambda^{-1} dx +c_{46}\int_{E\cap S_\lambda }g(c_{45}|x|^{1-N})|x|  dx
 \\
   &\le & c_1g (c_{45}\lambda^{N-2}) |E|+c_1m(\lambda )g (c_{45}\lambda^{N-2}) +c_1\int_{ \lambda }^\infty  m(s)dg (c_{45}s^{N-2})
   \\&&+c_{46}g (c_{45}\lambda^{N-1})\lambda^{-1} |E|+c_{46} m(\lambda )g (c_{45}\lambda^{N-1})\lambda^{-1}
   \\&& +c_{46}\int_{ \lambda }^\infty  m(s)d(g (c_{45}s^{N-1})s^{-1}).
\end{eqnarray*}
Since the critical index in (\ref{1.2}) is $\frac{N+1}{N-1}$, we have that
$$\int_{ \lambda }^\infty   g(c_{45}s^{N-2})d  m(s)=\lim_{T\to\infty}\int_{ \lambda }^{ T  }  g(js)d  m(s).
$$
Thus,
\begin{eqnarray*}
&&  m (\lambda) g(c_{45}\lambda^{N-2}) + \int_{ \lambda }^{ T  }    m(s)d  g(c_{45}s^{N-2})
 \\&=&  c_{47}T^{-N}g(c_{45}T^{N-2})+c_{48}\int_{ \lambda }^{ T  } s^{-1-N} g(c_{45}s^{N-2})ds
 \\&=&  c_{47}(T^{N-2})^{-\frac{N}{N-2}}g(c_{45}T^{N-2})+c_{49}\int_{(c_{45} \lambda)^{1/(N-2)} }^{ (c_{45} T)^{1/(N-2)}   } t^{-1-\frac{N}{N-2}} g(t)dt
\end{eqnarray*}
and
\begin{eqnarray*}
&&  m (\lambda) g(c_{45}\lambda^{N-1})\lambda^{-1}  + \int_{ \lambda }^\infty  m(s)d(g (c_{45}s^{N-1})s^{-1})
 \\&=&  c_{50}T^{-N-1}g(c_{45}T^{N-1})+c_{50}\int_{ \lambda }^{ T  } s^{-2-N} g(c_{45}s^{N-1})ds
 \\&=&  c_{50}(T^{N-1})^{-\frac{N+1}{N-1}}g(c_{45}T^{N-1})+c_{51}\int_{(c_{45} \lambda)^{1/(N-1)} }^{ (c_{45} T)^{1/(N-1)}   } t^{-1-\frac{N+1}{N-1}} g(t)dt,
\end{eqnarray*}
where $$ c_{47}(T^{N-2})^{-\frac{N}{N-2}}g(c_{45}T^{N-2})\to 0\quad{\rm as}\quad T\to\infty$$
and
$$c_{50}(T^{N-1})^{-\frac{N+1}{N-1}}g(c_{45}T^{N-1})\to 0\quad{\rm as}\quad T\to\infty$$
by the assumption (\ref{1.2}) and Lemma \ref{lm 4.1}.

Therefore,
$$ c_1m(\lambda )g (c_{45}\lambda^{N-2}) +c_1\int_{ \lambda }^\infty  m(s)dg (c_{45}s^{N-2}) \leq c_{49}\int_{(c_{45} \lambda)^{1/(N-2)} }^{ (c_{45} T)^{1/(N-2)}   } t^{-1-\frac{N}{N-2}} g(t)dt
$$
and
$$c_{46} m(\lambda )g (c_{45}\lambda^{N-1})\lambda^{-1} +c_{46}\int_{ \lambda }^\infty  m(s)d(g (c_{45}s^{N-1})s^{-1})\leq c_{51}\int_{(c_{45} \lambda)^{1/(N-1)} }^{ (c_{45} T)^{1/(N-1)}   } t^{-1-\frac{N+1}{N-1}} g(t)dt.$$
Notice that the quantities on the right-hand side tends to $0$
when $\lambda\to\infty$. The conclusion follows: for any
$\epsilon>0$, there exists $\lambda>0$ such that
$$c_{49}\int_{(c_{45} \lambda)^{1/(N-2)} }^{ (c_{45} T)^{1/(N-2)}   } t^{-1-\frac{N}{N-2}} g(t)dt \leq \frac{\epsilon}{6}
$$
and
$$c_{51}\int_{(c_{45} \lambda)^{1/(N-1)} }^{ (c_{45} T)^{1/(N-1)}   } t^{-1-\frac{N+1}{N-1}} g(t)dt \leq \frac{\epsilon}{6}.$$
For $\lambda$ fixed,  there exists $\delta>0$ such that
$$\int_E dx\leq \delta\Longrightarrow  c_1g (c_{45}\lambda^{N-2}) |E|\leq\frac{\epsilon}{6}\quad{\rm and}\quad c_{46}g (c_{45}\lambda^{N-1})\lambda^{-1} |E|\leq\frac{\epsilon}{6},
$$
 which implies that $\{\tilde g_n\circ \tilde v_n\}$ is uniformly integrable in
$L^1(B_1(0))$. Then
$$\tilde g_n\circ \tilde v_n \to g(w_{k}+\tilde v)-g(w_{k})\quad{\rm in}\quad L^1(B_1(0))$$
by Vitali convergence theorem.

Then passing to the limit as
$n\to +\infty$ in the identity (\ref{s52.1.1000}),
it implies that for any $\xi\in C^{1,1}_0(B_1(0))$,
$$
 \int_{B_1(0)} [w_k(-\Delta)\xi+g(w_k)\xi]dx=2k\frac{\partial \xi(0)}{\partial x_N}+j\xi(0).
$$
Thus, $v_{k,j}$ is a very weak solution of (\ref{eq 1.1}). The regularity results follows by Lemma \ref{lm 3.1}.

{\it Proof of $(ii)$. } It follows from (\ref{1.20}) and  (\ref{5.6}) that
\begin{eqnarray*}
 v_{k,j,n}(x)\le  2k\mathbb{G}_{B_1(0)}[\frac{\partial \delta_0}{\partial x_N}](x)+j \mathbb{G}_{B_1(0)}[ \delta_0 ](x)
 \end{eqnarray*}
 and
\begin{eqnarray*}
 v_{k,j,n}(x)  &\ge&  2k\mathbb{G}_{B_1(0)}[\frac{\partial \delta_0}{\partial x_N}](x)+j \mathbb{G}_{B_1(0)}[ \delta_0 ](x)-
 \mathbb{G}_{B_1(0)}[g_n(\mathbb{G}_{B_1(0)}[2k\frac{\partial \delta_0}{\partial x_N}+j\delta_0])](x)\\
     &\ge & 2k\mathbb{G}_{B_1(0)}[\frac{\partial \delta_0}{\partial x_N}](x)+j \mathbb{G}_{B_1(0)}[ \delta_0 ](x)-c_1w_g(x)-c_1  \mathbb{G}_{B_1(0)}[g(\mathbb{G}_{B_1(0)}[\delta_0])](x),
\end{eqnarray*}
where $w_g$  is the unique solution of (\ref{homo}) and
$$\lim_{|x|\to0^+}\frac{\mathbb{G}_{B_1(0)}[ \delta_0 ](x)}{\Gamma_N(x)}=1. $$
We see that  for any $e=(e_1,\cdots,e_N)\in \partial B_1(0)$ with $e_N\not=0$,
$$\lim_{t\to0^+}\mathbb{G}_{B_1(0)}[g(\mathbb{G}_{B_1(0)}[\delta_0])](te)t^{N-1}=0,$$
which, together with (\ref{3.001}), implies that
$$\lim_{t\to0^+}v_{k,j}(te)t^{N-1}=0.$$
 Then
$$
\lim_{t\to0^+}\mathbb{G}_{B_1(0)}[g(\mathbb{G}_{B_1(0)}[2k\frac{\partial  \delta_{0}}{\partial x_N}])](te)t^{N-1}=0.
$$
Thus,
$$\lim_{t\to0^+} \mathbb{G}_{B_1(0)}[\frac{\partial  \delta_{0}}{\partial x_N}](te)t^{N-1}=  e_N $$
 and by $x_N$-odd property of $w_{k}$,
we derive that
$$\lim_{t\to0^+} w_k(te)t^{N-1}=2k e_N. $$
 This ends the proof.\qquad$\Box$

\bigskip

\bigskip

\noindent{\bf Acknowledgements:}    H. Chen is supported by NSFC, No: 11401270, 11661045
and by the Jiangxi Provincial Natural Science Foundation, No: 20161ACB20007 and the
Project-sponsored by SRF for ROCS, SEM.

\end{document}